\documentclass [12pt,a4paper,reqno] {amsart}
\input amssymb.sty

\usepackage[dvips]{graphicx,color}
\usepackage{amsmath}
\usepackage{amsfonts}
\usepackage{amssymb}
\usepackage{bbm}
\usepackage[bbgreekl]{mathbbol}
\usepackage{latexsym}
\usepackage{epsfig}

\newtheorem{theorem}{Theorem}[section]
\newtheorem{lemma}[theorem]{Lemma}
\newtheorem{proposition}[theorem]{Proposition}
\newtheorem{corollary}[theorem]{Corollary}
\newtheorem{definition}[theorem]{Definition\rm}

\newtheorem{remark}{\it Remark\/}
\newtheorem{example}{\it Example\/}

\begin{document}

\newcommand{\fin}{$\Box$\\}
\newcommand{\ds}{\displaystyle}
\newcommand{\saut}[1]{\hfill\\[#1]}
\newcommand{\vsp}{\vspace{.15cm}}
\newcommand{\difrac}{\displaystyle \frac}
\newcommand{\dist}{\textrm{dist}}
\newcommand{\mbf}{\mathbf}
\newcommand{\sigmar}{\texttt{\large $\boldsymbol{\sigma}$}}
\newcommand{\alphar}{\texttt{\large $\boldsymbol{\alpha}$}}

\newcommand\alp{$\alpha_p$ }
\newcommand\all{$\alpha_l$ }
\newcommand\ho{H\"{o}lder }

\title[Stochastic $2$-microlocal analysis]{Stochastic $2$-microlocal analysis}
\author{Erick Herbin}
\address{Ecole Centrale Paris, Grande Voie des Vignes, 92295 Ch\^atenay-Malabry,
France} \email{erick.herbin@gmail.com}
\author{Jacques L\'evy-V\'ehel}
\address{Projet APIS, INRIA Saclay, Parc Orsay Universit\'e4 rue Jacques Monod - Bat P91893 Orsay Cedex, France}\email{jacques.levy-vehel@inria.fr}

\date{March 2008}

\subjclass[2000]{ 62\,G\,05, 60\,G\,15, 60\,G\,17, 60\,G\,18}

\begin{abstract}
A lot is known about the H\"older regularity of stochastic processes, in particular in the case of Gaussian processes.
Recently, a finer analysis of the local regularity of functions, termed {\it 2-microlocal analysis}, has been introduced in a deterministic frame: through the computation of the so-called {\it 2-microlocal frontier}, it allows in particular to predict the evolution of regularity under the action of (pseudo-) differential operators. 
In this work, we develop a 2-microlocal analysis for the study of certain stochastic processes. We show that moments of the increments allow, under fairly general conditions, to obtain almost sure lower bounds for the 2-microlocal frontier. 
In the case of Gaussian processes, more precise results may be obtained: 
the incremental covariance yields the almost sure value of the 2-microlocal frontier. As an application, we obtain new and refined regularity properties of fractional Brownian motion, multifractional Brownian motion, stochastic generalized Weierstrass functions, Wiener and stable integrals.
\end{abstract}

\maketitle

\section{Introduction}

It is a well-known fact that the pointwise \ho exponent of a function or
a process is not stable through the action of (pseudo-)differential operators.
In addition, it does not give a complete picture of the regularity of a function
at a given point, even if it is supplemented by the local \ho exponent (see
section \ref{s2ml} for definitions). Let us give two examples that show the significance of this for random processes.

\medskip

Define the function:
\begin{equation}\label{chirp}
c(t)=|t-t_0|^\gamma \sin(|t-t_0|^{-\beta}),
\end{equation}
where $t_0,\gamma$ and $\beta$ are positive real numbers.
This is an instance of a so-called ``chirp'', and is the simplest example of a function with non-trivial behaviour upon integro-differentiation (\cite{pspum}).
Now consider a multifractional Brownian motion $X_H(t)$ (see section \ref{gp2ml} for a definition). This process and its generalizations have been the subject of various 
studies in recent years (\cite{dozzi07,eh,peltier,ST1,MWX}). It is also currently used as a model in applications such as traffic engineering (\cite{mBmTCP}) or financial analysis (\cite{mBmFinance}).
It is parameterized by a function $H$, which controls its 
almost sure pointwise \ho exponent at each time. In applications, it is not expected that 
$H$ will behave smoothly. As a consequence $X_H(t)$ may have a complicated local 
regularity structure, and it is not easy to evaluate how it will be modified by various transforms of interest. 
Take for instance $H(t) = a + bc(t)$ in a neighbourhood of $t_0$. Here
$c$ is defined by (\ref{chirp}) and $a,b$ are chosen such that $c(t) \in (0,1)$ 
for $t$ in a neighbourhood of $t_0$.
Then the pointwise \ho exponent of $X_H(t)$ at $t_0$ is equal to $\gamma$ (provided that $H(t_0) > \gamma$). However, the pointwise \ho exponent of a fractional derivative of order $\varepsilon <\gamma/(1+\beta)$ of $X_H$ at $t_0$ is equal to $\gamma-\varepsilon(1+\beta)$. The fact that the regularity of $X_H$ decreases by more than $\varepsilon$ through $\varepsilon-$differentiation cannot be deduced from the sole knowledge of the exponent of $X_H$, but it is easily obtained with the help of the 2-microlocal analysis that we develop below. Example \ref{ExmBm} in section \ref{gp2ml} yields more details on this process, and describes a case where the evolution of the pointwise exponent is even ``stranger''.

As a second example, consider, for $t \geq 0$, the Wiener integral
\begin{equation}\label{EX}
Y_t=\int_0^t \eta(u) dB_u,
\end{equation}
where $B$ is standard Brownian motion and $\eta(t)=\sqrt{|c(t)|}$. The pointwise \ho exponent of $Y$ at $t_0$ is equal to $(\gamma + \beta + 1)/2$ and its local exponent at $t_0$ is equal to $\gamma/(2+2\beta) +1/2$, provided both these values are 
smaller than one. A fractional derivative of order $\varepsilon<\gamma/(2+2\beta) +1/2$ of $Y$ has pointwise and local exponents respectively equal to $(\gamma + \beta + 1)/2-\varepsilon(\beta + 1)$ and $\gamma/(2+2\beta) +1/2-\varepsilon$ at $t_0$. The variation of the pointwise exponent of $Y$ cannot be predicted from the sole knowledge of the exponents of $Y$, but are direct outputs of 2-microlocal analysis. See Example \ref{ExWi} in section \ref{gp2ml} for more on the Wiener integral (and example \ref{ExStable} which deals with the case of stable integrals).
\smallskip

Although these examples may seem somewhat {\it ad-hoc}, they allow to highlight, in a simple case, that there is more to local regularity that the mere local and pointwise \ho exponents classically considered. In addition, functions with even more irregular behaviour than the chirp can easily be exhibited, and such functions might well pop up in applications.
The aim of the present study is to provide new tools for the fine characterization of the regularity of stochastic processes, and in particular Gaussian processes, using the theory
of 2-microlocal analysis. Basically, 2-microlocal analysis allows to describe how the pointwise regularity of a function evolves under the action of (pseudo-) differential operators. This is useful in various areas, such as PDE (for which it was originally developed), signal or image analysis. The idea is to ``mix'', in a clever way, the local and pointwise characterizations of \ho regularity in a single condition involving two exponents (see inequality (\ref{spdf}) for a precise statement). So far, 2-microlocal analysis has only been considered in a deterministic frame. Since the pointwise regularity of random processes also is important both in theory and in applications such as the ones mentionned in the previous paragraph, it seems desirable to develop a stochastic version of 2-microlocal analysis.

We provide first steps in this direction below. We show first that an upper bound on moments on the increments of a process $X$ around a point $t \in \mathbf{R}^N_{+}$ provides an almost sure lower bound for the 2-microlocal frontier at this point. We also prove a related uniform result on $\mathbf{R}^N_{+}$. In the case where $X$ is a Gaussian process, we are able to obtain more precise results: 
the behaviour of the incremental covariance allows to obtain the almost sure value of the 2-microlocal frontier at ay given point. It also provides uniform results, which are however less precise. 
These results apply at once to classical processes, and allow for instance to recover easily known facts about multifractional Brownian motion. We obtain in addition new information which allow in particular to deal with the examples mentioned above concerning Wiener integrals and multifractional Brownian motion.

\medskip

The remaining of this paper is organized as follows: we start by recalling some basic facts about (deterministic) 2-microlocal analysis in section \ref{bf2ml}. Section \ref{sect1} contains our main results about the 2-microlocal analysis of continuous random processes: lower bounds for general processes, upper and lower bounds for Gaussian processes. We apply these results to various well-known processes in section \ref{gp2ml}. Finally, proofs of intermediate results are gathered in section \ref{boundsect}.

\section{Background: Deterministic $2$-microlocal analysis}\label{bf2ml}

2-microlocal analysis, which was introduced by J.M. Bony in \cite{Bony}, provides a tool that allows to predict the evolution of the local regularity of a function under the action of (pseudo-)differential operators. To be more precise, let $f^{(\varepsilon)}$ denote the fractional integral (when $\varepsilon <0$) or fractional derivative (when $\varepsilon >0$) of the real function $f$. The pointwise \ho exponent of $f^{(\varepsilon)}$ at $t$ is denoted $\alpha_{f^{(\varepsilon)}}(t)$ (see definition \ref{rkineq}). 
In several applications ({\it e.g.} PDE, signal or image processing), one needs to have access to the function $\mathcal{H}_t : \varepsilon\mapsto\alpha_{f^{(\varepsilon)}}(t)$.  Knowledge of $\mathcal{H}$ allows to answer questions such as: how much does one (locally) regularize the process $f$ by integrating it? The problem comes from the well-known fact that the pointwise \ho exponent is not stable under integro-differentiation: while it is true in simple situations that $\alpha_{f^{(n)}}(t)=\alpha_f(t) -n$, in general, one can only ensure that $\alpha_{f^{(n)}}(t) \leq \alpha_f(t) -n$. 2-microlocal analysis provides a way to assess the evolution of $\alpha_f$ through the use of a fine scale of functional spaces. These {\it 2-microlocal spaces}, denoted $C^{s,s'}$, generalize the classical \ho spaces in a way we describe now.

Since $\mathcal{H}$ cannot be deduced from the sole knowledge of $\alpha_f$\footnote{Because of the inequality $\alpha_{f^{(n)}}\leq\alpha_f -n$, $\mathcal{H}$ has to decrease faster than $\varepsilon\mapsto -\varepsilon$. One can show that, apart from this and a certain regularity property, there are no other constraints on the evolution of the pointwise \ho exponent (see \cite{pspum}).}, predicting changes in the regularity of a process under integro-differentiation basically requires recording the whole function $\varepsilon\mapsto\alpha_{f^{(\varepsilon)}}(t)$. 
2-microlocal analysis does this in a clever way and without having to compute any integro-differentials: 
it associates to any given point $t$ a curve in a abstract space, its {\it 2-microlocal frontier}, whose slope is the rate of increase of $\varepsilon \mapsto \alpha_{f^{(\varepsilon)}}(t)$. The 2-microlocal frontier may be estimated through a fine analysis of the local regularity of $f$ around $t$. This analysis can be conducted in the Fourier (\cite{Bony}), wavelet (\cite{Jaff}) or time (\cite{KKJLV,pspum}) domains. We shall use in this work the time-domain characterization of 2-microlocal spaces. Proofs of the results of this section and more information on 2-microlocal analysis may be found in \cite{Bony,pspum}.

\begin{definition}[Time domain definition of 2-microlocal spaces]\label{def2ml}
Let $x_0\in \mathbf{R}$, and $s,s'$ be two real
numbers satisfying $s+s'>0$, $s+s'\not \in \mathbf{N}$, and $s'<0$
(and thus $s > 0$). Let $m=[s+s']$ (the integer part of $s+s'$).

A function $f: \hspace{1mm} \mathbf{R}\rightarrow\mathbf{R}$
belongs to $C^{s,s'}_{x_0}$ if and only if its $m^{th}$ derivative
exists around $x_0$, and if there exist a positive real $\delta$, a
polynomial $P$ of degree not larger than $[s]-m$, and a constant $C$,
that verify
\begin{equation*}
\left | \frac{\partial^mf(x)-P(x)}{|x-x_0|^{[s]-m}}-
\frac{\partial^mf(y)-P(y)} {|y-x_0|^{[s]-m}} \right | \leq
C|x-y|^{s+s'-m}(|x-y|+|x-x_0|)^{-s'-[s]+m}
\end{equation*}
for all $x,y$ such that $0<|x-x_0| < \delta$, $0<|y-x_0| <
\delta$.
\end{definition}

Except in this introductory section, we shall restrict to the case where $(s,s')$ verify $0 < s+s' < 1$, $s < 1$, $s' < 0$. This corresponds to the situation where $f$ is not differentiable at $x_0$, but has some global regularity in the neighbourhood of $x_0$. More precisely, we shall assume that there exists an interval $I$ containing $x_0$ and a real number $\eta \in (0,1)$ such that $f$ belongs to the global \ho space $C^{\eta}(I)$. 
This restriction allows to avoid certain technicalities in the analysis. We believe all the results should hold in the general case with appropriate modifications. When $(s,s')$ satisfy the above inequalities, $m$ in definition \ref{def2ml} is equal to 0 and the polynomial $P$ is a constant. 
As a consequence, the inequality characterizing 2-microlocal spaces reduces to:

\begin{definition}[Time domain definition of 2-microlocal spaces, case of non-differentiable continuous functions]\label{def2mlsimple}
Let $D = \{(s,s') : 0 < s+s' < 1, s < 1, s' < 0 \}$.  A
function $f: \hspace{1mm} \mathbf{R}\rightarrow\mathbf{R}$ belongs to $C^{s,s'}_{x_0}$, with $(s,s') \in D$, if there exist a positive real $\delta$ and a constant $C$ such that for all $(x,y)$ with $0<|x-x_0|<\delta$, $0<|y-x_0| < \delta$,
\begin{equation}\label{spdf}
|f(x)-f(y)|\leq C|x-y|^{\sigma}(|x-x_0|+|x-y|)^{-s'},
\end{equation}
where $\sigma=s+s'$.
\end{definition}

Recall that the pointwise \ho exponent of $f$ at $x$ is defined as the supremum of the $\alpha$ such that $f$ belongs to pointwise \ho spaces $C^{\alpha}_x$. 2-microlocal spaces use two parameters $(s,s')$. The relevant notion generalizing the pointwise exponent is the {\it 2-microlocal frontier}.
In order to define this frontier, consider the {\it 2-microlocal domain} of $f$ at $x_0$, i.e. the set $E(f,x_0)=\{(s,s') : f\in C^{s,s'}_{x_0}\}$. One can prove that $E(f,x_0)$ is always a convex subset of the abstract plane $(s,s')$. 
The 2-microlocal frontier $\Gamma(f,x_0)$ is the convex curve in the $(s,s')$-plane defined by
\begin{align*}
 \Gamma(f,x_0) :\quad & \mathbf{R} \to \mathbf{R}\\
& s'\mapsto  s(s')=\sup\left\{r: f\in C^{r,s'}_{x_0} \right\}.
\end{align*}

For various reasons (see \cite{pspum}), it is useful to describe the 2-microlocal frontier as a function $s' \mapsto \sigma(s')$ (recall that $\sigma = s+s'$), and this is the parameterization we shall mainly use in the following. By abuse of language, we shall refer to $\sigma(s')$ as the 2-microlocal frontier in the sequel.\\

The following property of $\sigma(s')$ will be useful:
\begin{proposition}
\label{frontdef} The 2-microlocal frontier of $f$ at any point $x_0$, seen
as a function $s' \mapsto \sigma(s')$, verifies
\begin{itemize}
\item
$\sigma(s')$ is a concave, non-decreasing function,
\item
$\sigma(s')$ has left and right derivatives always between 0 and
1.
\end{itemize}
\end{proposition}

We present now the fundamental properties of the 2-microlocal frontier (see Definition \ref{rkineq} for the definition of the local \ho exponent $\tilde{\alpha}_f(x_0)$ of $f$ at $x_0$):

\begin{proposition}[Stability under fractional integro-differentiation]\label{stab}
For any function $f:\mathbf{R}\rightarrow\mathbf{R}$,
for all $ (s,s')\in\mathbf{R}$, for all $x_0$ and for all $\varepsilon$
\begin{equation*}
f\in C^{s,s'}_{x_0} \Longleftrightarrow  f^{(\varepsilon)} \in C^{s-\varepsilon,s'}_{x_0}.
\end{equation*}
\end{proposition}

\begin{proposition}[Pointwise \ho exponent]\label{pexp}
Assume $f\in C^{\eta}(\mathbf{R})$ for some $\eta>0$.\\
Then, the pointwise \ho exponent of $f$ at $x_0$ is given by 
\begin{equation*}
\alpha_f(x_0)=- \inf\left\{s': \sigma(s')\geq 0 \right\},
\end{equation*} 
with the convention that $\alpha_f(x_0)=+\infty$ if $\sigma(s')\geq 0$ for
all $s'$.
\end{proposition}

\begin{proposition} [Local \ho exponent]\label{lexp}
Assume $f\in C^{\eta}(\mathbf{R})$ for some $\eta>0$.
Then, the local \ho exponent of $f$ at $x_0$ is given by $$\tilde{\alpha}_f(x_0)= \sigma(0).$$
\end{proposition}

\begin{figure}
\label{frontier}
\begin{center}
\includegraphics[width=9cm,height = 9cm]{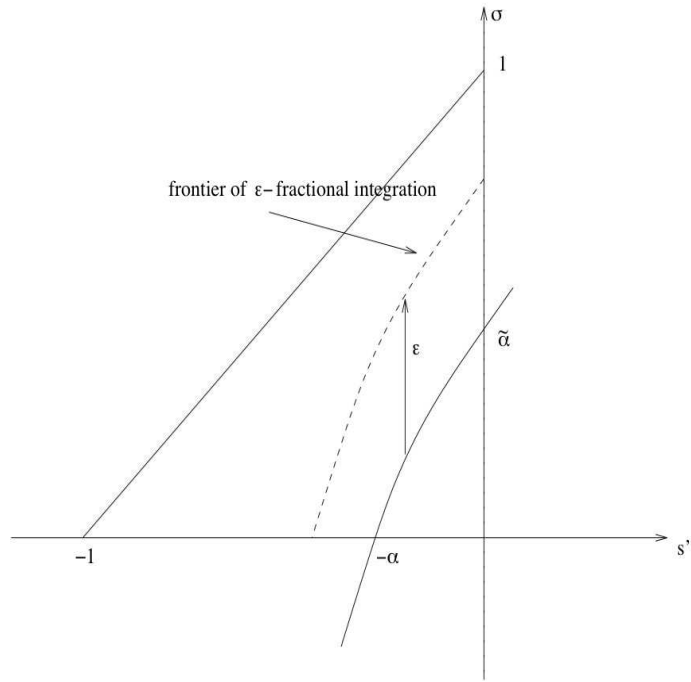}
\caption{2-microlocal frontier of a function $f$ in the
$(s',\sigma)$-plane (solid curve). The intersection with the $s'$ axis
occurs at $-\alpha$, and the one with the $\sigma$ axis occurs at $\tilde \alpha$. 
The dotted curve shows the frontier of an $\varepsilon-$ integral of $f$.
It is obtained by translation of $\varepsilon$ along the ordinate axis.} 
\end{center}
\end{figure}

The above propositions show that {\it the 2-microlocal frontier contains the whole information pertaining to $\varepsilon\mapsto\alpha_{f^{(\varepsilon)}}(t)$}. Indeed, fractional integro-differentiation of order $\varepsilon$ amounts to translating the 2-microlocal frontier by $-\varepsilon$ along the $\sigma$ direction in the $(s',\sigma)$ plane (proposition \ref{stab}). The pointwise \ho exponent of $f^{(\varepsilon)}$  is then given by (minus) the intersection of the translated frontier with the $s'$ axis (proposition \ref{lexp}), provided $\varepsilon >\tilde{\alpha}_f$ (proposition \ref{pexp}). See figure \ref{frontier}.

For our examples below, we will need the 2-microlocal frontiers of the
following functions: the function $x\mapsto |t-t_0|^{\gamma}$ has a trivial frontier at $t_0$: 
it is parallel to the first bisector and passes through the point $(-\gamma,0)$.
The 2-microlocal frontier of the chirp $t \mapsto |t-t_0|^{\gamma} \sin(\frac{1} {|t-t_0|^{\beta}})$ at $t_0$ is the straight line defined by $\sigma(s') = \frac{1} {\beta+1} s' +\frac{\gamma}{\beta+1}$. Finally, the Weierstrass function $W_s(t)=\sum_{n=1}^{+\infty}\lambda^{-nh}\sin(\lambda^n t)$, where $\lambda \geq 2$, $0<h<1$ has the same frontier at all $t$: 
$\sigma(s')=s'+h$ for $s'\leq 0$, and $\sigma(s)'=h$ for $s'\geq 0$. See \cite{pspum}.

This ends our recalls on 2-microlocal analysis.

\section{$2$-microlocal analysis of random processes}\label{sect1}

In the remaining of this article, we shall perform the 2-microlocal analysis of certain random processes. We start by transposing the notion described in the previous section in a stochastic frame, and by defining some quantities that will prove useful for computing the almost-sure frontier of our processes.

\subsection{Stochastic $2$-microlocal analysis}\label{s2ml}

Let $X=\left\{X_t;\;t\in\mathbf{R}^N_{+}\right\}$ be a 
continuous random process.
For each $t_0\in\mathbf{R}^N_{+}$, let us define the
{\it $2$-microlocal frontier} of $X$ at $t_0$
as the random function $s'\mapsto\sigmar_{t_0}(s')$, defined
for  $s'\in(-\infty;0)$ by
\begin{equation}\label{ladef}
\sigmar_{t_0}(s')=\sup\left\{\sigma;\; \limsup_{\rho\rightarrow 0}
\sup_{t,u\in B(t_0,\rho)}\frac{|X_t-X_u|}{\|t-u\|^{\sigma}\rho^{-s'}}
<\infty \right\}.
\end{equation}
Each couple $\left(s';\sigmar_{t_0}(s')\right)$ could be called a
``$2$-microlocal exponent'' of $X$ at $t_0$.

While (\ref{ladef}) makes sense for all $s'\in(-\infty;0)$ and all $\sigma$, $\sigmar$ will {\it not} in general coincide with a stochastic version of the 2-microlocal frontier: indeed, it states that, for a given realization $\omega$, $X(\omega) \in C^{s,s'}_{t_0}$ whenever $s < \sigmar_{t_0}(s') - s'$ {\it only when definition \ref{def2mlsimple} may be applied in place of the more general definition \ref{def2ml}}. Thus, in the sequel, we shall always assume that $\sigmar_{t_0}$ (or, more correctly, its representation in the $(s,s')$ plane) intersects the region $D = \{(s,s') : 0 < s+s' < 1, s < 1, s' < 0 \}$. 
Recall that this is equivalent to assuming that $X$ is not differentiable at $t_0$, but belongs to $C^{\eta}(I)$ for some interval $I$ containing $t_0$ and some $\eta \in (0,1)$. 
Note also that $D$ may equivalently by characterized by
\begin{equation*}
D = \left\{(s',\sigma) : -1<\hspace{0.05cm}s'<\hspace{0.05cm}0,
\hspace{0.4cm} 0< \sigma < 1+s' \right\}.
\end{equation*}
This assumption will allow us to avoid technicalities entailed by the use of definition \ref{def2ml}.
However, the fact that we restrict to this region implies that we will only be able to predict the variations of the pointwise exponent through {\it differentiation} of order $\varepsilon$ {\it small enough}. Indeed, integration means translating the frontier in the direction of positive $\sigma$. As a consequence, evaluating the pointwise exponent after integration requires knowledge of at least a part of the frontier that is below the $s'$ axis prior to integration. Furthermore, all the results are valid only when the local exponent remains positive, and this means that we cannot consider differentiation of too high an order.

\smallskip
\noindent
In order to obtain almost sure results,
we shall seek triples $(\eta,\mu,\nu)$ such that
\begin{equation}\label{Cgendet}
\limsup_{\rho\rightarrow 0}\sup_{t,u\in B(t_0,\rho)}
\frac{E\left[X_t-X_u\right]^\eta}{\|t-u\|^{N+\mu}\rho^{-\nu}}<+\infty.
\end{equation}

In the case where $X$ is Gaussian, it is natural to take $\eta=2$ and to consider the set
of couples $\left(s',\sigma\right)$ such that
\begin{equation}\label{Cdet}
\limsup_{\rho\rightarrow 0}\sup_{t,u\in B(t_0,\rho)}
\frac{E\left[X_t-X_u\right]^2}{\|t-u\|^{2\sigma}\rho^{-2s'}}<+\infty.
\end{equation}

This approach is usual in the analysis of Gaussian processes (see, e.g., \cite{adler,dudley,talagrand96}), and it leads naturally to define ``deterministic 2-microlocal spaces'' as follows:

\begin{definition}
A Gaussian process $X$ is said to belong to $\mathbb{C}^{s,s'}_{t_0}$ for
a fixed $t_0\in\mathbf{R}^N_{+}$ and some $s,s'$ such that
\begin{equation}\label{detsp}
\left\{
\begin{array}{c}0 < s+s' < 1\\ s < 1 \\ s' < 0 \end{array}
\right.
\end{equation}
if  condition (\ref{Cdet}) is satisfied for $\sigma=s+s'$.
\end{definition}

Recall the definitions of the pointwise and local H\"older exponents of $X$ at $t_0$:
\begin{definition}\label{rkineq}
The pointwise and local \ho exponents of $X$ at $t_0$ are defined as the random variables:
\begin{align*}
\alphar(t_0)=\sup\left\{\alpha; \limsup_{\rho\rightarrow 0}
\sup_{t,u\in B(t_0,\rho)}\frac{|X_t-X_u|}{\rho^{\alpha}}
<\infty \right\}, \\
\widetilde{\alphar}(t_0)=\sup\left\{\alpha; \limsup_{\rho\rightarrow 0}
\sup_{t,u\in B(t_0,\rho)}\frac{|X_t-X_u|}{\|t-u\|^{\alpha}}
<\infty \right\}.
\end{align*}
Note that although these quantities are in general random variables, we will omit the dependency in $\omega$ and write $\alphar(t_0)$ and $\widetilde{\alphar}(t_0)$ instead of $\alphar(t_0,\omega)$ and $\widetilde{\alphar}(t_0,\omega)$. 
\end{definition}
It is easily seen that for all $t_0$ and all $s'<0$, we have (for all $\omega$)
\begin{equation}\label{locfrpoint}
\widetilde{\alphar}(t_0)\leq\sigmar_{t_0}(s')-s'\leq\alphar(t_0).
\end{equation}

To show this inequality, proceed as follows:
\begin{itemize}
\item For all $0<\sigma<\sigmar_{t_0}(s')$,
\begin{eqnarray*}
\frac{|X_t-X_u|}{\rho^{\sigma-s'}}
= \frac{|X_t-X_u|}{\|t-u\|^{\sigma}\rho^{-s'}}
\left(\frac{\|t-u\|}{\rho}\right)^{\sigma}
\end{eqnarray*}
which gives $\sigma-s'\leq\alphar(t_0)$. As a consequence, $\sigmar_{t_0}(s')-s'\leq\alphar(t_0)$.
\item For all $\alpha<\widetilde{\alphar}(t_0)$,
\begin{eqnarray*}
\frac{|X_t-X_u|}{\|t-u\|^{\alpha+s'}\rho^{-s'}}
=\frac{|X_t-X_u|}{\|t-u\|^{\alpha}}
\left(\frac{\|t-u\|}{\rho}\right)^{-s'}
\end{eqnarray*}
which gives $\alpha+s'\leq\sigmar_{t_0}(s')$. As a consequence,
$\widetilde{\alphar}(t_0)\leq\sigmar_{t_0}(s')-s'$.
\end{itemize}
\fin

Section \ref{sectlowerbound} is concerned with general processes and obtains in this case lower bounds for the $2$-microlocal frontier. In section \ref{sectupperbound}, we focus on Gaussian processes and obtain upper bounds. We would like to mention here that
there is a large body of results about the regularity of sample paths of Gaussian processes. A non-exhaustive list of authors that contributed to these results includes Berman, Dudley, Fernique, Marcus, Orey, Rosen, Pitt, Pruitt, Talagrand, Xiao...
We refer to \cite{adler, davar, ledouxtalagrand, marcusrosen} for a contemporary and complete survey of these works. This field of research is still very active (e.g. \cite{cras, dozzi03, xiao}).
Extremely precise results are known, under various conditions, on uniform moduli of continuity and local moduli of continuity, using the terminology of \cite{marcusrosen}, chapter 7.
Uniform moduli of continuity give information which is finer than the mere local \ho exponent, and $\widetilde{\alphar}$ is easily obtained once an exact uniform modulus of continuity is known. Likewise, local moduli of continuity give richer information than the pointwise \ho exponent, and $\alphar$ may be deduced at once from an exact local modulus of continuity. Theorem 7.1.2 in \cite{marcusrosen} is a general, rather abstract, result giving a uniform modulus of continuity based on a majorizing measure (see also \cite{adler, ledouxtalagrand, talagrand96}). 
Precise bounds for moduli of continuity are given in \cite[Theorem 7.2.1]{marcusrosen} under an assumption related to, but weaker than, stationarity of the increments. Theorem 7.4.5 gives another result under different but fairly general conditions. When the Gaussian process has stationary increments and when its incremental variance verifies some regularity assumptions, exact moduli of continuity may be obtained, see Theorems 7.2.9, 7.2.10, 7.2.14, 7.2.15, 7.6.4 and 7.6.9 in \cite{marcusrosen}.
These results do not, however, apply directly to our situation, for two reasons.

First, for the processes we will mostly be interested in, it is not straightforward to obtain a majorizing measure, so it is not clear how to apply Theorem 7.1.2 of \cite{marcusrosen}.
As for the other Theorems mentioned above, the ones using stationarity of the increments cannot be used since the processes we have in mind are ``strongly'' non increment-stationary: mBm, for instance, is such that its increments of any order are never stationary as soon as the $H$ function is not constant.
The other, weaker, assumptions made in these theorems are either that $E\left[X_{t}-X_{u}\right]^2 \leq \varphi(|t-u|)$ for some strictly increasing function $\varphi$, or that $E\left[X_{t}X_{u}\right] \leq \min(E\left[X_{t}\right]^2, E\left[X_{u}\right]^2)$. None of these assumptions is verified by mBm.

Second, although the results of \cite{marcusrosen} mentioned above, when available, give much more precise information on the pointwise and local regularity than the ones we are going to obtain, they would not directly allow to compute the 2-microlocal frontier. As a consequence, they would be of no use for, {\it e.g.}, computing the variation of regularity under differentiation of the ``chirped multifractional Brownian motion'' mentioned at the beginning of the introduction. Indeed, as explained at the beginning of this work, the sole knowledge of $\widetilde{\alphar}$ and $\alphar$, or even of exact moduli of continuity, is not sufficient in order to predict the evolution of pointwise regularity under integro-differentiation, and one needs rather to obtain bounds ``mixing'' the two exponents, as in inequality (\ref{spdf}).
Although it should be possible to generalize the approach described in \cite{marcusrosen} based on the powerful tool of majorizing measures to obtain 2-microlocal characterizations, we follow below a different path. In our frame, the Gaussian assumption is not necessary to obtain lower bounds, and we deal instead with the more general class of processes satisfying condition (\ref{chara}).

\subsection{Lower bound for the $2$-microlocal frontier of stochastic processes}\label{sectlowerbound}

In this section, we give conditions for the paths of a stochastic process, not necessarily Gaussian, to belong to a given $2$-microlocal space $C^{s,s'}_{t_0}$.

\subsubsection{Pointwise almost sure result}

\begin{proposition} \label{propmin}
Let $X=\left\{X_t;\;t\in\mathbf{R}_{+}^N\right\}$ be a cadlag (right continuous with left limits) stochastic process.
Assume that for some $t_0$, there exists two constants $C>0$ and $\rho_0>0$
such that
\begin{equation}\label{chara}
\forall 0<\rho<\rho_0, \forall t,u\in B(t_0,\rho);\quad
E\left|X_t-X_u\right|^{\eta} \leq C\|t-u\|^{N+\mu} \rho^{-\nu}
\end{equation}
with $\eta>0$, $\mu>0$ and $\nu<0$.\\
Then, almost surely, the paths of the process $X$ belong to
$C^{\sigma-\nu/\eta,\nu/\eta}_{t_0}$, for all $\sigma<\frac{\mu}{\eta}$.

In other words, the $2$-microlocal frontier of $X$ at $t_0$ verifies
\begin{equation*}
\sigmar_{t_0}\left(\frac{\nu}{\eta}\right)\geq \frac{\mu}{\eta} \quad\textrm{a.s.}
\end{equation*}
\end{proposition}

\proof
Let $s'=\nu/\eta$ and $0<\sigma<\mu/\eta$.
Let us take $\rho=2^{-n}$ for $n\geq n_0=-\log_2\rho_0$ and set
$D^m_n(t_0)=\left\{t_0+k.2^{-(m+n)};k\in\left\{0, \pm 1,\dots,\pm 2^m\right\}^N\right\}$.
Let us consider the event
\[
\Omega^m_n=\left\{\max_{\scriptstyle k,l\in\left\{0,\dots,\pm 2^m\right\}^N
\atop \scriptstyle \|k-l\|=1}
\left|X_{t_0+k.2^{-(m+n)}}-X_{t_0+l.2^{-(m+n)}}\right|>
2^{-\sigma (m+n)}2^{s'n}\right\}.
\]
We have
\begin{align*}
P\left\{\Omega^m_n\right\}
&\leq\sum_{\scriptstyle k,l\in\left\{0,\dots,\pm 2^m\right\}^N
\atop \scriptstyle \|k-l\|=1}P\left\{
\left|X_{t_0+k.2^{-(m+n)}}-X_{t_0+l.2^{-(m+n)}}\right|>
2^{-\sigma (m+n)}2^{s'n}\right\}\\
&\leq\sum_{\scriptstyle k,l\in\left\{0,\dots,\pm 2^m\right\}^N
\atop \scriptstyle \|k-l\|=1}\frac{E\left|X_{t_0+k.2^{-(m+n)}}-X_{t_0+l.2^{-(m+n)}}\right|^{\eta}}
{2^{-\eta\sigma (m+n)}2^{\eta s'n}}\\
&\leq 2N\;C\; \underbrace{\#\left\{0,\dots,\pm 2^m\right\}^N}_{\left(1+2^{(m+1)}\right)^N}
2^{-(N+\mu-\eta\sigma) (m+n)}.
\end{align*}
Then,
\begin{equation*}
P\left\{\Omega^m_n\right\}
\leq (2N.2^{2N})\;C\; 2^{-(\mu-\eta\sigma)m}\; 2^{-(N+\mu-\eta\sigma) n}
\end{equation*}
and
\begin{align*}
P\left\{\exists m;\Omega^m_n\right\}=P\left\{\bigcup_m \Omega^m_n\right\}
&\leq\sum_m P\left\{\Omega^m_n\right\}\\
&\leq\frac{(2N.2^{2N})\;C\; 2^{-(N+\mu-\eta\sigma) n}}{1-2^{-(\mu-\eta\sigma)}}.
\end{align*}
The Borel-Cantelli lemma implies existence of a random variable
$n^{*}\geq n_0$ such that, almost surely,
\[
\forall n\geq n^{*}, \forall m\in\mathbf{N};\quad
\max_{\scriptstyle k,l\in\left\{0,\dots,\pm 2^m\right\}^N
\atop \scriptstyle \|k-l\|=1}
\left|X_{t_0+k.2^{-(m+n)}}-X_{t_0+l.2^{-(m+n)}}\right|
\leq 2^{-\sigma (m+n)}2^{s'n}.
\]
Therefore, by induction, we get for all $n\geq n^{*}$ and all $m\in\mathbf{N}$
\begin{align*}
\forall q>m,\; \forall t,u\in D_n^q(t_0);\textrm{ s.t. } \|t-u\|<2^{-(m+n)};&\\
|X_t-X_u|&\leq
2\left(\sum_{j=m+1}^{q}2^{-\sigma (j+n)}\right) 2^{s'n}\\
&\leq\frac{2.2^{-\sigma (m+n+1)}}{1-2^{-\sigma}}
2^{s' n}
\end{align*}
which leads to
\begin{align*}
\forall t,u\in D_n(t_0)=\bigcup_m D^m_n(t_0);\quad 
|X_t-X_u|\leq\frac{2}{1-2^{-\sigma}} \|t-u\|^{\sigma}\ 2^{s'n}
\end{align*}
and using the right continuity of $X$,
\begin{align*}
\forall t,u\in B(t_0,2^{-n});\quad 
|X_t-X_u|\leq\frac{2}{1-2^{-\sigma}} \|t-u\|^{\sigma}\ 2^{s'n}.
\end{align*}
Hence, almost surely, for all $\rho\in (0,2^{-n^{*}})$, there exists
$n>n^{*}$ such that $2^{-(n+1)}\leq\rho\leq 2^{-n}$ and
\begin{align}
\forall t,u\in B(t_0,\rho); \quad
|X_t-X_u|&\leq
\frac{2}{1-2^{-\sigma}}
\|t-s\|^{\sigma}\ 2^{s'n}\nonumber\\
&\leq 2^{-s'}\frac{2}{1-2^{-\sigma}}
\|t-u\|^{\sigma} \rho^{-s'}.
\end{align}
This inequality proves that the sample paths of $X$ belong to $C^{\sigma-\nu/\eta,\nu/\eta}_{t_0}$.
This implies that for all $\epsilon>0$, $P\{\sigmar_{t_0}(\nu/\eta)\geq\mu/\eta-\epsilon\}=1$. Taking $\epsilon\in\mathbf{Q}_+$, we get 
$P\{\sigmar_{t_0}(\nu/\eta)\geq\mu/\eta\}=1$.
\fin

\begin{corollary} \label{cormingauss}
Let $X=\left\{X_t;t\in\mathbf{R}_{+}^N\right\}$ be a continuous Gaussian process.
Assume that for some $t_0$, there exists two constants $C>0$ and $\rho_0>0$
such that
\begin{equation}\label{charagauss}
\forall 0<\rho<\rho_0, \forall t,u\in B(t_0,\rho);\quad
E\left[X_t-X_u\right]^2\leq C\|t-u\|^{2\sigma} \rho^{-2s'}
\end{equation}
with $\sigma>0$ and $s'<0$.\\
Then, almost surely, the paths of the process $X$ belong to
$C^{\tilde{\sigma}-s',s'}_{t_0}$, for all $\tilde{\sigma}<\sigma$.\\
In other words, 
\begin{align*}
X\in\mathbb{C}^{\sigma-s',s'}_{t_0}
\Rightarrow \left[ \forall \tilde{\sigma}<\sigma,\quad
X\in C^{\tilde{\sigma}-s',s'}_{t_0} \quad\textrm{a.s.} \right]
\end{align*}
and consequently,
\begin{equation*}
\sigmar_{t_0}(s')\geq\sigma \quad\textrm{a.s.}
\end{equation*}
\end{corollary}

\proof
From (\ref{charagauss}), for all $p\in\mathbf{N}^*$ we get
\begin{equation*}
\forall 0<\rho<\rho_0, \forall t,u\in B(t_0,\rho);\quad
E\left[X_t-X_u\right]^{2p}\leq C\; \lambda_p\;\|t-u\|^{2p\sigma} \rho^{-2ps'}
\end{equation*}
where $\lambda_p$ is the positive constant such that for all centered Gaussian random variable $Y$, we have $E[Y^{2p}]=\lambda_p \left(E[Y^2]\right)^p$.

For any $p\in\mathbf{N}^*$ such that $2p\sigma > N$, proposition \ref{propmin} implies that the sample paths of $X$ belong to $C_{t_0}^{\tilde{\sigma}-s', s'}$ for all $\tilde{\sigma}<\sigma-\frac{N}{2p}$. The result follows.

\fin

\begin{remark}
The assumptions of corollary \ref{cormingauss} are equivalent to the existence of
$\sigma>0$ and $s'<0$ such that
\[
\limsup_{\rho\rightarrow 0}\sup_{t,u\in B(t_0,\rho)}
\frac{E\left[X_t-X_u\right]^2}{\|t-u\|^{2\sigma}\rho^{-2s'}}<+\infty.
\]
\end{remark}

In proposition \ref{propmin}, if (\ref{chara}) holds for all $t_0$ in $[a,b]$, by Kolmogorov's criterion, the process $X$ admits a version $\tilde{X}$, which is continuous on $[a,b]$. Therefore, the cadlag condition may be thought unnecessary:
The regularity of the paths of $X$ would stand for the regularity of the paths of one of its continuous version.\\
On the contrary to pointwise results, uniform almost sure regularity at $t_0\in\mathbf{R}^N_+$ should depend on the considered version. In that view, we assume that $X$ is cadlag and 
as all the cadlag versions of $X$ are indistinguishable, the question of uniformity in $t_0\in\mathbf{R}^N_+$ is worth further investigation.

\subsubsection{Uniform almost sure result on $\mathbf{R}^N_{+}$}

Inequalities (\ref{locfrpoint}) show that an uniform lower bound for the
local H\"older exponent gives an uniform lower bound for
$\sigmar_{t_0}(s')$.

\begin{proposition}\label{propminunif}
Let $X=\left\{X_t;\;t\in\mathbf{R}^N_{+}\right\}$ be a continuous stochastic
process. Assume that there exist a constant $\eta>0$ and a positive function $\mu$ which admits a positive lower bound, and that for all $t_0\in\mathbf{R}^N_{+}$, there exists
$C_0>0$ and $\rho_0>0$ such that
\begin{equation}\label{eqhypth}
\forall t,u\in B(t_0,\rho_0),
E\left|X_t-X_u\right|^{\eta}\leq C_0\; \|t-u\|^{N+\mu(t_0)}
\end{equation}
Then, almost surely
\begin{equation}
\forall t_0\in\mathbf{R}^N_{+};\quad\widetilde{\alphar}(t_0)\geq
\liminf_{u\rightarrow t_0}\frac{\mu(u)}{\eta}
\end{equation}
and, as a consequence
\begin{equation}
\forall t_0\in\mathbf{R}_{+}^N,\forall s'<0;\quad\sigmar_{t_0}(s')\geq
s'+\liminf_{u\rightarrow t_0}\frac{\mu(u)}{\eta}.
\end{equation}
\end{proposition}

\proof
First, let us suppose that the function $\mu$ is constant.\\
By (\ref{eqhypth}) and Kolmogorov's criterion, for all
$t_0\in\mathbf{R}_{+}^N$, there exists a modification
$\tilde{X}_{t_0}$ of $X$ which is $\alpha$-H\"older continuous for all
$\alpha\in (0,\mu/\eta)$ on $B(t_0,\rho_0)$.
Therefore the local H\"older exponent of $\tilde{X}_{t_0}$ satisfy
\[
\forall t\in B(t_0,\rho_0);\quad\widetilde{\alphar}_{\tilde{X}_{t_0}}(t)\geq\frac{\mu}{\eta}.
\]
As a consequence, for all $t_0\in\mathbf{R}_{+}^N$, there exists $\rho_0>0$ such
that
\[
P\left\{\forall t\in B(t_0,\rho_0);\;
\widetilde{\alphar}(t)\geq\frac{\mu}{\eta}\right\}=1.
\]
For all $a,b\in\mathbf{Q}_{+}^N$, such that $a\prec b$, we have
\[
[a,b]\subset\bigcup_{t_0\in [a,b]}B(t_0,\rho_0).
\]
As $[a,b]$ is compact, there exists a finite number of balls $B_1,\dots,B_n$
such that
\[
[a,b]\subset\bigcup_{i=1}^n B_i
\]
and
\[
\forall i=1,\dots,n;\quad P\left\{\forall t\in
B_i;\;\widetilde{\alphar}(t)\geq\frac{\mu}{\eta}\right\}=1.
\]
Therefore, we get
\[
P\left\{\forall t\in [a,b];\;
\widetilde{\alphar}(t)\geq\frac{\mu}{\eta}\right\}=1.
\]
As $\mathbf{R}_{+}^N$ can be covered by a countable number of compact sets
$[a,b]$, this leads to
\begin{equation}\label{minunifconst}
P\left\{\forall t\in\mathbf{R}_{+}^N;\;
\widetilde{\alphar}(t)\geq\frac{\mu}{\eta}\right\}=1.
\end{equation}
Using (\ref{locfrpoint}) and the continuity of $s'\mapsto\sigmar_{t_0}(s')$,
the result follows.
\vsp

In the general case of a non-constant function $\mu$,
for all $a,b\in\mathbf{Q}_{+}^N$
with $a\prec b$ and all $\sigma=\inf_{u\in [a,b]}\mu(u)/\eta-\epsilon$ with
$\epsilon>0$,
there exists a constant $C>0$ such that
\[
\forall t,u\in [a,b];\quad E\left|X_t-X_u\right|^{\eta}\leq C\; \|t-u\|^{N+\sigma}.
\]
Then (\ref{minunifconst}) implies the existence of
$\Omega^{*}\subset\Omega$ such that $P\left\{\Omega^{*}\right\}=1$ and for all
$\omega\in\Omega^{*}$,
\[
\forall a,b\in\mathbf{Q}_{+}^N,\forall\epsilon\in\mathbf{Q}_{+},
\forall t_0\in\overbrace{[a,b]}^{\circ};\quad
\widetilde{\alphar}(t_0)\geq\inf_{u\in [a,b]}\frac{\mu(u)}{\eta}-\epsilon.
\]
Therefore, taking two sequences $\left(a_n\right)_{n\in\mathbf{N}}$ and
$\left(b_n\right)_{n\in\mathbf{N}}$ such that $\forall
n\in\mathbf{N};\;a_n<t_0<b_n$ and converging to $t_0$,
we have for all $\omega\in\Omega^{*}$
\begin{equation*}
\forall t_0\in\mathbf{R}^N_{+};\quad\widetilde{\alphar}(t_0)\geq
\liminf_{u\rightarrow t_0}\frac{\mu(u)}{\eta}.
\end{equation*}
\fin

{\it Example : Linear Multifractional Stable Motion on $\mathbf{R}$

The Linear Multifractional Stable Motion (LMSM) is an $\alpha-$stable process
(with $0 < \alpha < 2$) which is obtained from the better-known linear
fractional stable motion by replacing the constant $H$ by a H\"older-continuous function
$H(t)$. See \cite[Chapter 3]{ST} for details on the linear
fractional stable motion and \cite{ST1,ST2} for the definition and properties of
the LMSM. Let us consider such an LMSM, with the following assumptions : $\alpha >1$,
and, for all $t$, $1/\alpha < H(t) < 1$. Assume in addition that, for all $t_0$,
$\tilde \alpha_H(t_0) > 1/\alpha$. 
It is shown in \cite{ST2}, inequality (3.6), that:
\begin{equation*}\label{lmsm}
\forall p \in (1,\alpha), \forall t,u\in B(t_0,\rho_0),
E\left|X_t-X_u\right|^p\leq C(p) 
|t-u|^{p(\tilde \alpha_H(t_0)\wedge \min_{v \in B(t_0,\rho_0)} H(v))}.
\end{equation*}
A direct application of Proposition \ref{propminunif} then yields:
\begin{equation*}
\forall t_0\in\mathbf{R}_{+}^N,\forall s'<0;\quad\sigmar_{t_0}(s')\geq
s'+ \tilde \alpha_H(t_0)\wedge H(t_0) - \frac {1}{\alpha}.
\end{equation*}
\cite{ST2} contains more precise statements on the local and pointwise \ho
exponents of the LMSM. See section \ref{intstoch} for related results on stable integrals.
}
\medskip

Using the same technique as in Corollary \ref{cormingauss}, one gets, in the
particular case of Gaussian processes:

\begin{corollary}\label{corminunifgauss}
Let $X=\left\{X_t;\;t\in\mathbf{R}^N_{+}\right\}$ be a continuous Gaussian
process. Assume that there exists a function $\sigma$ which admits a positive
lower bound, and that for all $t_0\in\mathbf{R}^N_{+}$, there exists
$C_0>0$ and $\rho_0>0$ such that
\begin{equation}\label{eqhypth}
\forall t,u\in B(t_0,\rho_0),
E\left[X_t-X_u\right]^2\leq C_0\; \|t-u\|^{2\sigma(t_0)}.
\end{equation}
Then, almost surely
\begin{equation}
\forall t_0\in\mathbf{R}^N_{+};\quad\widetilde{\alphar}(t_0)\geq
\liminf_{u\rightarrow t_0}\sigma(u)
\end{equation}
and, as a consequence
\begin{equation}
\forall t_0\in\mathbf{R}_{+}^N,\forall s'<0;\quad\sigmar_{t_0}(s')\geq
s'+\liminf_{u\rightarrow t_0}\sigma(u).
\end{equation}
\end{corollary}

The uniform almost sure lower bound for the $2$-microlocal frontier can be
improved when $E\left|X_t-X_u\right|^{\eta}$ admits, in the ball
$B(t_0,\rho)$, an upper bound function of $\|t-u\|$ and $\rho$,
uniformly in $t_0$. We only state the result in the Gaussian case.

\begin{proposition}\label{propunifminsigma}
Let $X=\left\{X_t;\;t\in\mathbf{R}^N_{+}\right\}$ be a continuous Gaussian
process.
Assume that for all $t_0$, there exists a real function $\varsigma_{t_0}$
such that
\[
\forall t_0, \forall s'<0;\quad \varsigma_{t_0}(s')\geq m(s')>0
\]
and that there exists two constants $C>0$ and $\delta>0$ such that
for all $t_0 \in [0,1]$ and all $s'<0$
\begin{equation}
\forall 0<\rho<\delta,\;\forall t,u\in B(t_0,\rho);\quad
E\left[X_t-X_u\right]^2\leq C \; \|t-u\|^{2\varsigma_{t_0}(s')} \rho^{-2s'}.
\end{equation}
Then for all $s'<0$, almost surely
\[
\forall t_0\in [0,1];\quad
\sigmar_{t_0}(s')\geq\liminf_{u\rightarrow t_0}\varsigma_{u}(s').
\]
\end{proposition}

The proof of Proposition \ref{propunifminsigma} is an improvement of the one of
Proposition \ref{propmin} and is somewhat technical. It is given in section \ref{boundsect}.

\subsection{Upper bound for the $2$-microlocal frontier of Gaussian processes}\label{sectupperbound}

In this section, we assume that the considered processes are Gaussian. 
We show that an almost sure upper bound can be established for 
the $2$-microlocal frontier of the sample paths.

\subsubsection{Pointwise almost sure result}
To get an almost sure upper bound for $\sigmar_{t_0}(s')$, we need the following

\begin{lemma} \label{majlem}
Let $X=\left\{X_t;\;t\in\mathbf{R}^N_{+}\right\}$ be a continuous Gaussian
process.
Assume that for some $t_0\in\mathbf{R}^N_{+}$, there exist $\sigma>0$ and $s'<0$
such that there are two
sequences $\left(t_n\right)_{n\in\mathbf{N}}$ and
$\left(u_n\right)_{n\in\mathbf{N}}$ converging to $t_0$, and a constant $c>0$
such that
\begin{equation*}
\forall n\in\mathbf{N};\quad E\left[X_{t_n}-X_{u_n}\right]^2
\geq c \|t_n-u_n\|^{2\sigma}\rho_n^{-2s'}
\end{equation*}
where $\rho_n\geq\max(\|t_n-t_0\|,\|u_n-t_0\|)$. Then the $2$-microlocal
exponent satisfies almost surely
\[
\sigmar_{t_0}(s') \leq \sigma.
\]
\end{lemma}

\proof
Let $\epsilon>0$ and consider two sequences $(t_n)_{n\in\mathbf{N}}$
and $(u_n)_{n\in\mathbf{N}}$ as in the statement of the Lemma.\\
For all $n\in\mathbf{N}$, the law of the random variable
$\frac{X_{t_n}-X_{u_n}}{\|t_n-u_n\|^{\sigma+\epsilon}\rho_n^{-s'}}$ is
$\mathcal{N}(0,\sigma_n^2)$.\\
>From the assumption, we have $\sigma_n\geq c\|t_n-u_n\|^{-2\epsilon}\rightarrow +\infty$
as $n\rightarrow +\infty$.\\
Then, for all $\lambda>0$,
\begin{align*}
P\left\{\frac{\|t_n-u_n\|^{\sigma+\epsilon}\rho_n^{-s'}}{|X_{t_n}-X_{u_n}|}<\lambda
\right\}&=P\left\{\frac{|X_{t_n}-X_{u_n}|}{\|t_n-u_n\|^{\sigma+\epsilon}\rho_n^{-s'}}
>\frac{1}{\lambda}\right\}\\
&=\int_{|x|>\frac{1}{\lambda}}\frac{1}{\sqrt{2\pi}\sigma_n}
\exp\left(-\frac{x^2}{2\sigma_n^2}\right).dx\\
&=\frac{1}{2\pi}\int_{|x|>\frac{1}{\lambda\sigma_n}}
\exp\left(-\frac{x^2}{2}\right).dx\stackrel{n\rightarrow +\infty}
{\longrightarrow} 1.
\end{align*}
Therefore the sequence
$\left(\frac{\|t_n-u_n\|^{\sigma+\epsilon}\rho_n^{-s'}}{|X_{t_n}-X_{u_n}|}\right)_{n\in\mathbf{N}}$
converge to $0$ in probability. This implies the
existence of a subsequence which converge
to $0$ almost surely.
Thus, almost surely, $\sigmar_{t_0}(s')\leq\sigma+\epsilon$. Taking
$\epsilon\in\mathbf{Q}_{+}$, the result follows.
\fin

\begin{remark}\label{remupp}
The assumptions of lemma \ref{majlem} are equivalent to the existence of
$\sigma>0$ and $s'<0$ such that
\begin{equation*}
\limsup_{\rho\rightarrow 0}\sup_{t,u\in B(t_0,\rho)}
\frac{E\left[X_t-X_u\right]^2}{\|t-u\|^{2\sigma}\rho^{-2s'}}
>0.
\end{equation*}
\end{remark}

\subsubsection{Uniform almost sure result on $\mathbf{R}^N_{+}$}

In contrast with what happens in the case of the lower bound, 
if the assumptions of lemma \ref{majlem} are satisfied for all $t_0$, 
the conclusion holds uniformly in $t_0$:

\begin{proposition}\label{propmajunif}
Let $X=\left\{X_t;t\in\mathbf{R}^N_{+}\right\}$ be a continuous Gaussian process.
Suppose that the assumptions of lemma \ref{majlem} are satisfied for all
$t_0\in\mathbf{R}_{+}^N$,
with the same $\sigma$ and $s'$.
Then the $2$-microlocal exponent satisfies almost surely
\[
\forall t_0\in\mathbf{R}^N_{+};\quad \sigmar_{t_0}(s')\leq\sigma.
\]
\end{proposition}

For the sake of readability, we postpone the proof of 
this Proposition to  section \ref{boundsect}. 
In the light of Remark \ref{remupp}, Proposition \ref{propmajunif} implies

\begin{theorem}\label{thmajunif}
Let $X=\left\{X_t;\;t\in\mathbf{R}^N_{+}\right\}$ be a continuous Gaussian
process. Assume that for all $t_0\in\mathbf{R}^N_{+}$ and all $s'<0$, there
exists $\varsigma_{t_0}(s')>0$ such that
\[
\limsup_{\rho\rightarrow 0}\sup_{t,u\in B(t_0,\rho)}
\frac{E\left[X_t-X_u\right]^2}
{\|t-u\|^{2\varsigma_{t_0}(s')}\rho^{-2s'}}
>0.
\]
Then for all $s'<0$, we have almost surely
\[
\forall t_0\in\mathbf{R}^N_{+};\quad
\sigmar_{t_0}(s')\leq\limsup_{u\rightarrow t_0}\varsigma_{u}(s').
\]
\end{theorem}

\proof
Under the assumptions of the theorem, for all $a,b\in\mathbf{Q}^N_{+}$
with $a\prec b$ and all
$\varsigma(s')=\sup_{u\in [a,b]}\varsigma_{u}(s')$,
\[
\forall t_0\in[a,b];\quad\limsup_{\rho\rightarrow 0}\sup_{s,t\in B(t_0,\rho)}
\frac{E\left[X_t-X_s\right]^2}{\|t-s\|^{2\varsigma(s')}\rho^{-2s'}}>0.
\]

Therefore, by proposition \ref{propmajunif}, for all $s'<0$, there exists a set
$\Omega^{*}\subset\Omega$ with
$P\left\{\Omega^{*}\right\}=1$ such that for all $\omega\in\Omega^{*}$
\[
\forall a,b\in\mathbf{Q}^N_{+},
\forall t_0\in\overbrace{[a,b]}^{\circ};\quad
\sigmar_{t_0}(s')\leq\sup_{u\in [a,b]}\varsigma_{u}(s').
\]
Taking two sequences $\left(a_n\right)_{n\in\mathbf{N}}$ and
$\left(b_n\right)_{n\in\mathbf{N}}$ such that $\forall
n\in\mathbf{N};\;a_n<t_0<b_n$ and converging to $t_0$,
we have for all $\omega\in\Omega^{*}$
\begin{equation*}
\forall t_0\in\mathbf{R}^N_{+};\quad\sigmar_{t_0}(s')\leq
\limsup_{u\rightarrow t_0}\varsigma_{u}(s').
\end{equation*}
\fin

\subsection{Where is, almost surely, the $2$-microlocal frontier?}
In this section, we show that, not surprisingly, the $2$-microlocal frontier of
the paths of a Gaussian process can be evaluated by studying its incremental
covariance.
The proofs rely on the computation of almost sure lower and upper bounds for the
frontier, which were developed in sections \ref{sectlowerbound} and
\ref{sectupperbound}.

As a counterpart to the random H\"older exponents and frontier
$\sigmar_{t_0}(s')$, let us introduce the {\it deterministic
local H\"older exponent}
\begin{equation}
\widetilde{\mathbb{\bbalpha}}(t_0)=
\sup\left\{\alpha;\;\limsup_{\rho\rightarrow 0}\sup_{t,u\in B(t_0,\rho)}
\frac{E\left[X_t-X_u\right]^2}{\|t-u\|^{2\alpha}}<\infty
\right\}
\end{equation}
and the {\it deterministic $2$-microlocal frontier}
$s'\mapsto\mathbb{\bbsigma}_{t_0}(s')$:
\begin{align}\label{sigdet}
\mathbb{\bbsigma}_{t_0}(s')&=
\sup\left\{\sigma;\;\limsup_{\rho\rightarrow 0}\sup_{t,u\in B(t_0,\rho)}
\frac{E\left[X_t-X_u\right]^2}{\|t-u\|^{2\sigma}\rho^{-2s'}}<\infty
\right\} \\
&=\sup\left\{\sigma;\; X \in \mathbb{C}^{\sigma-s',s'}_{t_0}\right\}. \nonumber
\end{align}

The same proof as in the frame of deterministic functions allows
to show that the deterministic 2-microlocal frontier
$s'\mapsto\mathbb{\bbsigma}_{t_0}(s')$ is concave and thus
continuous on $(-1,0)$.

The two following sections give almost sure pointwise and uniform
results concerning the 2-microlocal frontier (recall that we
assume that $\sigmar_{t_0}$ intersects the region defined by
conditions (\ref{detsp})).

\subsubsection{Pointwise almost sure $2$-microlocal frontier}
Proposition \ref{propmin} in section \ref{sectlowerbound} shows that for
all $s'$, almost surely,
$$\mathbb{\bbsigma}_{t_0}(s')\leq\sigmar_{t_0}(s').$$  
Conversely, according to lemma \ref{majlem}, for all $s'$, almost surely,
$$\sigmar_{t_0}(s')\leq\mathbb{\bbsigma}_{t_0}(s').$$ 
Using additionally the continuity of the frontier, one may thus state

\begin{theorem}\label{thone}
Let $X=\left\{X_t;\;t\in\mathbf{R}^N_{+}\right\}$ be a continuous Gaussian
process.
For any $t_0\in\mathbf{R}^N_{+}$, the $2$-microlocal frontier of $X$ at $t_0$
is almost surely equal to the function $s'\mapsto\mathbb{\bbsigma}_{t_0}(s')$.
\end{theorem}

In the same way as in the deterministic frame, the almost sure values of the pointwise and local H\"older exponents can be computed from the almost sure $2$-microlocal frontier (see Propositions \ref{pexp} and \ref{lexp}).

\begin{corollary}\label{pe}
For any $t_0\in\mathbf{R}^N_{+}$, the pointwise \ho exponent of $X$ at
$t_0$ is almost surely equal to $- \inf \{s': \mathbb{\bbsigma}_{t_0}(s')\geq 0\}$,
provided $\mathbb{\bbsigma}_{t_0}(0)>0$.
\end{corollary}

\begin{corollary}\label{le}
For any $t_0\in\mathbf{R}^N_{+}$, the local \ho exponent of $X$ at
$t_0$ is almost surely equal to $\mathbb{\bbsigma}_{t_0}(0)$, provided
$\mathbb{\bbsigma}_{t_0}(0)>0$.
\end{corollary}

\subsubsection{Uniform almost sure result on $\mathbf{R}^N_{+}$}\label{parfrunif}
Proposition \ref{propminunif} and theorem \ref{thmajunif} in sections \ref{sectlowerbound} and \ref{sectupperbound} provide some almost sure results
about the $2$-microlocal frontier and the local H\"older exponent,
uniformly in $t_0\in\mathbf{R}_{+}^N$.

\begin{theorem}\label{aslocal}
Let $X=\left\{X_t;\;t\in\mathbf{R}^N_{+}\right\}$ be a continuous Gaussian
process.
Assume that the function
$t_0\mapsto\liminf_{u\rightarrow t_0}\tilde{\mathbb{\bbalpha}}(u)$ is positive.
Then, almost surely
\begin{equation}
\forall t_0\in\mathbf{R}_{+}^N;\quad
\liminf_{u\rightarrow t_0}\widetilde{\mathbb{\bbalpha}}(u) \leq \widetilde{\alphar}(t_0) \leq
\limsup_{u\rightarrow t_0}\widetilde{\mathbb{\bbalpha}}(u).
\end{equation}
\end{theorem}

\proof
By definition of $\widetilde{\mathbb{\bbalpha}}(t_0)$, for all $\epsilon>0$, and all
$t_0\in\mathbf{R}^N_{+}$,
there exist $C_0>0$ and $\rho_0>0$ such that
\[
\forall t,u\in B(t_0,\rho_0);\quad
E\left[X_t-X_u\right]^2\leq C_0 \|t-u\|^{2(\tilde{\mathbb{\bbalpha}}(t_0)-\epsilon)}.
\]
Proposition \ref{propminunif} implies that, almost surely,
\begin{equation*}
\forall t_0\in\mathbf{R}^N_{+};\quad
\widetilde{\alphar}(t_0)\geq
\liminf_{u\rightarrow t_0}\widetilde{\mathbb{\bbalpha}}(u) - \epsilon
\end{equation*}
and, taking $\epsilon\in\mathbf{Q}_{+}$,
\begin{equation*}
\forall t_0\in\mathbf{R}^N_{+};\quad
\widetilde{\alphar}(t_0)\geq
\liminf_{u\rightarrow t_0}\widetilde{\mathbb{\bbalpha}}(u).
\end{equation*}

Conversely,
using theorem \ref{thmajunif} with $s'=0$, for all $\epsilon>0$,
we have almost surely
\begin{equation*}
\forall t_0\in\mathbf{R}^N_{+};\quad
\widetilde{\alphar}(t_0)\leq
\limsup_{u\rightarrow t_0}\tilde{\mathbb{\bbalpha}}(u) + \epsilon.
\end{equation*}
Taking $\epsilon\in\mathbf{Q}_{+}$, we get
\begin{equation*}
\forall t_0\in\mathbf{R}^N_{+};\quad
\widetilde{\alphar}(t_0)\leq
\limsup_{u\rightarrow t_0}\tilde{\mathbb{\bbalpha}}(u).
\end{equation*}
\fin

\begin{corollary}
Let $X=\left\{X_t;\;t\in\mathbf{R}^N_{+}\right\}$ be a continuous Gaussian
process.
Assume that the function
$t_0\mapsto\widetilde{\mathbb{\bbalpha}}(t_0)$ is continuous and positive.
Then, almost surely
\begin{equation}
\forall t_0\in\mathbf{R}_{+}^N;\quad
\widetilde{\alphar}(t_0)=\tilde{\mathbb{\bbalpha}}(t_0).
\end{equation}
\end{corollary}

By remark \ref{rkineq}, theorems \ref{thmajunif} and \ref{aslocal} imply

\begin{corollary} \label{corfrontunif}
Let $X=\left\{X_t;\;t\in\mathbf{R}^N_{+}\right\}$ be a Gaussian process.
Assume that the function
$t_0\mapsto\widetilde{\mathbb{\bbalpha}}(t_0)$ is continuous and positive.
Then, almost surely
\begin{equation}
\forall t_0\in\mathbf{R}^N_{+},\forall s'<0;\quad
\widetilde{\mathbb{\bbalpha}}(t_0)+s'\leq\sigmar_{t_0}(s')\leq
\limsup_{u\rightarrow t_0}\mathbb{\bbsigma}_{u}(s').
\end{equation}
\end{corollary}

Theorem \ref{aslocal} and corollary \ref{corfrontunif} only give
bounds for the uniform almost sure $2$-microlocal frontier.
However, in specific cases, we are able to obtain a sharp result:
This happens for instance for fractional Brownian motion and
regular multifractional Brownian motion, as we show below.

\section{Applications to some well-known Gaussian and non-Gaussian processes}\label{gp2ml}

The results of the previous section can be used to compute the
almost sure 2-microlocal frontier of some well-known Gaussian
processes such as (multi)fractional Brownian motion,
generalized Weierstrass function, and Wiener integrals. We also briefly
consider the case of stable integrals in section \ref{intstoch}.

\subsection{Fractional Brownian motion}

Fractional Brownian motion (fBm) is one of the simplest processes whose regularity has been deeply studied (see \cite{adler} for a recent account). It is defined as the continuous Gaussian process $B^H=\left\{B^H_t;\;t\in\mathbf{R}_{+}\right\}$ such that for all $t,u\in\mathbf{R}_{+}$,
\begin{equation}\label{eqfbm}
E\left[B^H_t-B^H_u\right]^2=|t-u|^{2H}
\end{equation}
where $H\in (0,1)$.\\
It is well-known that fBm is almost surely H\"{o}lder-continuous
but nowhere differentiable. As a consequence, its $2$-microlocal
frontier intersects the region defined by conditions
(\ref{detsp}). The results of paragraph \ref{parfrunif} can then
be applied to fBm. Theorem \ref{aslocal} directly yields the value
of the almost sure local H\"older exponent uniformly on
$\mathbf{R}_{+}$. The uniformity of (\ref{eqfbm}) in the whole of
$\mathbf{R}_{+}$, then allows to get the almost sure
$2$-microlocal frontier of fBm, uniformly in $\mathbf{R}_{+}$.

\begin{proposition} \label{propfbm}
Almost surely,
the 2-microlocal frontier at any $t_0$ of the fractional Brownian motion in the region
\[
\left\{
\begin{array}{c}  0 < \sigma < 1 + s' \\ -1 <  s' < 0 \end{array}
\right.
\]
is equal to the line $\sigma=H+s'$.
\end{proposition}

\proof
According to Corollary \ref{corfrontunif}, the result relies on the fact that $\mathbb{\bbsigma}_{t_0}(s')= H+s'$ for all $t_0\in\mathbf{R}_{+}$ and all $s'<0$. 

To prove this fact, we first observe that definition of fBm implies that
\begin{equation*}\label{fbmmin}
\forall t_0\in\mathbf{R}_{+};\quad
\widetilde{\mathbb{\bbalpha}}(t_0)=H,
\end{equation*}
which gives the lower bound for $\mathbb{\bbsigma}_{t_0}(s')$.
Secondly, to get the upper bound, we consider any sequences
$\left(\rho_n\right)_{n\in\mathbf{N}}$ converging to $0$, $\left(u_n\right)_{n\in\mathbf{N}}$ and $\left(t_n\right)_{n\in\mathbf{N}}$ such that for all $n\in\mathbf{N}$, $u_n,t_n\in B(t_0,\rho_n)$ and $|t_n-u_n|=\rho_n$. 
As a consequence,
for all $n$,
\begin{align*}
\frac{E\left[X_{t_n}-X_{u_n}\right]^2}{|t_n-u_n|^{2(H+s')}\rho_n^{-2s'}}
&=|t_n-u_n|^{-2s'} \rho_n^{2s'}\\
&=1
\end{align*}
which gives
\[
\limsup_{\rho\rightarrow 0}\sup_{t,u\in B(t_0,\rho)}
\frac{E\left[X_t-X_u\right]^2}{|t-u|^{2(H+s')}\rho^{-2s'}} > 0
\]
and then,
\begin{equation*}\label{fbmmaj}
\forall t_0\in\mathbf{R}_{+},\forall s'<0;\quad
\mathbb{\bbsigma}_{t_0}(s')\leq H+s'.
\end{equation*}
\fin

As the pointwise (resp. local) H\"older exponent is the intersection of the 2-microlocal
frontier with the axis $s'=0$ (resp. the line $\sigma=0$),
one recover the following well-known results as immediate consequences of Proposition \ref{propfbm}.
\begin{corollary}\label{corfBm}
The local and pointwise H\"older exponents satisfy almost surely
\[
\forall t_0\in\mathbf{R}_{+};\quad
\widetilde{\alphar}(t_0)=\alphar(t_0)=H.
\]
\end{corollary}

\subsection{Multifractional Brownian motion}
As shown in corollary \ref{corfBm}, the local regularity of fBm is
constant along the paths. A natural extension of fBm is to substitute the
constant parameter $H$, with a function $t\mapsto H(t)$ taking values
in $(0,1)$. This leads to
multifractional Brownian motion (see \cite{benassi}, \cite{peltier}).
The mBm can be defined as the process $X=\left\{X_t; t\in\mathbf{R}_{+}\right\}$
such that
\[
X_t=\int_{\mathbf{R}}\left[|t-u|^{H(t)+\frac{1}{2}}-|u|^{H(t)+\frac{1}{2}}
\right].\mathbb{W}(du)
\]
where $\mathbb{W}$ denotes the white noise of $\mathbf{R}$.\\

This process is now well studied, with many results on its local regularity (\cite{peltier}), higher dimensional versions (\cite{eh,MWX}), extension to the case where $H$ is itself random (\cite{AT}), study of its local time (\cite{dozzi03,dozzi07,MWX}) and more.
The covariance structure of mBm was first computed in \cite{cohencov} and the asymptotic behavior of the incremental covariance 
was investigated in \cite{eh}.
For all $a,b\in [0,1]$, and all $t_0\in [a,b]$, there exist positive constants
$K(t_0)$ and $L(t_0)$ such that
\begin{align}\label{approx}
\forall t,u\in B(t_0,\rho);\quad
E\left[X_t-X_u\right]^2 = K(t_0)\|t-u\|^{2H(t)}
+L(t_0)\left[H(t)-H(u)\right]^2\nonumber\\
+ o_{a,b}\left(\|t-u\|^{2H(t)}\right)
+ o_{a,b}\left(H(t)-H(u)\right)^2.
\end{align}
This approximation, together with the fact that mBm is H\"{o}lder-continuous but not differentiable, allows to compute the almost sure 2-microlocal frontier of the mBm at any point $t_0$. However, in contrast to the case of fBm, (\ref{approx}) only gives local information about the covariance. As a consequence, obtaining almost sure results uniformly in $t_0$ requires further work.

\subsubsection{Pointwise almost sure $2$-microlocal frontier of the mBm}

\begin{proposition} \label{propmbm}
The 2-microlocal frontier of the multifractional Brownian
motion in the region
\[
\left\{
\begin{array}{c} -1 < s' < 0 \\ 0 < \sigma < 1 + s' \end{array}
\right.
\]
is, at any fixed $t_0$, almost surely equal to the ``minimum'' of the
2-microlocal frontier of $H$ at $t_0$ and the line
$\sigma=H(t_0)+s'$. More precisely, for all $t_0$,
$\sigmar_{t_0}(s')=(H(t_0)+s') \wedge \beta_{t_0}(s')$ with
probability one, where $\beta_{t_0}(s')$ denotes the 2-microlocal
frontier of the deterministic function $H$ at $t_0$.
\end{proposition}

\proof
By definition, for each $s'\in (-\infty; 0)$,
\begin{equation}
\beta_{t_0}(s')=\sup\left\{\beta;\; \limsup_{\rho\rightarrow 0}
\sup_{t,u\in B(t_0,\rho)}\frac{|H(t)-H(u)|}{\|t-u\|^{\beta}\rho^{-s'}}
<\infty \right\}.
\end{equation}
We have to distinguish the 2 following cases:
\begin{itemize}
\item $H(t_0)+s'<\beta_{t_0}(s')$\\
For all $\sigma<H(t_0)+s'$, there exists $\eta_0>0$ such that
\[
\forall t\in B(t_0,\eta_0);\quad\sigma<H(t)+s'.
\]
Then, for all $0<\rho<\eta_0$
\begin{align*}
\frac{\|t-u\|^{2H(t)}}{\|t-u\|^{2\sigma}\rho^{-2s'}}=
\frac{\|t-u\|^{2(H(t)-\sigma)}}{\rho^{-2s'}}\leq
\frac{(2\rho)^{2(H(t)-\sigma)}}{\rho^{-2s'}}
\rightarrow 0
\end{align*}
and
\begin{align*}
\frac{\left[H(t)-H(u)\right]^2}{\|t-u\|^{2\sigma}\rho^{-2s'}}\rightarrow 0.
\end{align*}
Then by (\ref{approx}), we have $\sigma\leq\mathbb{\bbsigma}_{t_0}(s')$.
This implies
\begin{equation}\label{mbmcas1min}
H(t_0)+s'\leq\mathbb{\bbsigma}_{t_0}(s').
\end{equation}
For all $\sigma$ s.t. $H(t_0)+s'<\sigma<\beta_{t_0}(s')$,
there exists $\eta_1>0$ such that
\[
\forall t\in B(t_0,\eta_1);\;H(t)+s'<\sigma.
\]
Let us consider
$\rho_n=\frac{1}{n}$ and $u_n,t_n\in B(t_0,\rho_n)$ such that
$|t_n-u_n|=\rho_n$. For $\frac{1}{n}\leq\eta_1$, we have
\begin{align*}
\frac{\|t_n-u_n\|^{2H(t_n)}}{\|t_n-u_n\|^{2\sigma}\rho_n^{-2s'}}=
\rho_n^{2(H(t_n)+s'-\sigma)}\rightarrow +\infty
\end{align*}
and
\begin{align*}
\frac{\left[H(t_n)-H(u_n)\right]^2}{\|t_n-u_n\|^{2\sigma}\rho_n^{-2s'}}\rightarrow 0.
\end{align*}
Then from (\ref{approx}), we get $\sigma\geq\mathbb{\bbsigma}_{t_0}(s')$.
This implies
\begin{equation}\label{mbmcas1maj}
\mathbb{\bbsigma}_{t_0}(s')\leq H(t_0)+s'
\end{equation}

\item $\beta_{t_0}(s')<H(t_0)+s'$\\
There exists $\eta_2>0$ such that
\[
\forall t\in B(t_0,\eta_2);\quad\beta_{t_0}(s')<H(t)+s'.
\]
For all $\sigma<\beta_{t_0}(s')$, we have
\begin{align*}
\frac{\|t-u\|^{2H(t)}}{\|t-u\|^{2\sigma}\rho^{-2s'}}\rightarrow 0
\end{align*}
and
\begin{align*}
\frac{\left[H(t)-H(u)\right]^2}{\|t-u\|^{2\sigma}\rho^{-2s'}}\rightarrow 0.
\end{align*}
Then, by (\ref{approx}), we have $\sigma\leq\mathbb{\bbsigma}_{t_0}(s')$.
This implies
\begin{equation}\label{mbmcas2min}
\beta_{t_0}(s')\leq\mathbb{\bbsigma}_{t_0}(s').
\end{equation}
For all $\sigma$ s.t. $\beta_{t_0}(s')<\sigma<H(t_0)+s'$,
there exist sequences:
\begin{itemize}
\item $(\rho_n)_n$ of positive real numbers converging to $0$,
\item $(u_n)_n$ and $(t_n)_n$ s.t. $\forall n; u_n,t_n\in B(t_0,\rho_n)$
\end{itemize}
such that
\begin{align*}
\frac{\left[H(t_n)-H(u_n)\right]^2}{\|t_n-u_n\|^{2\sigma}\rho_n^{-2s'}}
\rightarrow +\infty.
\end{align*}
Moreover, there exists $\eta_3>0$ such that
\[
\forall t\in B(t_0,\eta_3);\quad\sigma<H(t)+s'.
\]
Let $N_3\in\mathbf{N}$ be such that $\forall n\geq N_0;\;0<\rho_n<\eta_3$.
For $n\geq N_3$, we have
\begin{align*}
\frac{\|t_n-u_n\|^{2H(t_n)}}{\|t_n-u_n\|^{2\sigma}\rho_n^{-2s'}}\rightarrow 0.
\end{align*}
Then, from (\ref{approx}), we get $\sigma\geq\mathbb{\bbsigma}_{t_0}(s')$.
This implies
\begin{equation}\label{mbmcas2maj}
\mathbb{\bbsigma}_{t_0}(s')\leq\beta_{t_0}(s').
\end{equation}
\end{itemize}
From (\ref{mbmcas1min}), (\ref{mbmcas1maj}), (\ref{mbmcas2min}) and
(\ref{mbmcas2maj}), the result follows with theorem \ref{thone}.
\fin

As in fBm's case, the almost sure pointwise and local H\"older exponents of mBm can be deduced from its almost sure $2$-microlocal frontier.

\begin{corollary}
At any $t_0$, the pointwise and local \ho exponents of the
multifractional Brownian motion verify almost surely:
$$\alphar(t_0)=H(t_0) \wedge \beta(t_0)$$
$$\widetilde{\alphar}(t_0)=H(t_0) \wedge \tilde{\beta}(t_0)$$
where $\beta(t_0)$ and $\tilde{\beta}(t_0)$
denote the pointwise and local \ho exponents of $H$ at $t_0$.
\end{corollary}
This result was already stated in \cite{eh}.

\begin{example} \label{ExmBm}
Following the introduction section, consider an mBm $X$ whose $H$ function is equal to $a + b|t-t_0|^\gamma \sin(|t-t_0|^{-\delta})$ in a neighbourhood of $t_0$, where $\gamma>0, \delta>0$, and $a,b$ are chosen such that $H(t) \in (0,1)$ for $t$ in a neighbourhood of $t_0$.
We may apply proposition \ref{propmbm} to get the almost sure 2-microlocal frontier of $X$ at $t_0$. 

Recall from section \ref{bf2ml} that the 2-microlocal frontier of $H$ is the function $\beta_{t_0}:s'\mapsto\beta_{t_0}(s')=\frac{1} {\delta+1} s' +\frac{\gamma}{\delta+1}$. Recall also that, because we are restricting to the part of the frontier that intersects the region defined by (\ref{detsp}), we can only consider what happens through differentiation of small enough order. 

Three situations may occur, depending on the values of $a$, $\gamma$ and $\delta$. If $a=H(t_0)<\frac{\gamma}{\delta+1}$, then the 2-microlocal frontier of $X$ is the line $s'\mapsto a+s'$, and nothing interesting happens from the 2-microlocal point of view.
If $a > \gamma$, then the frontier of $X$ is equal to $\beta_{t_0}$.
As a consequence, differentiation of $X$ of order $\varepsilon$ decreases its pointwise exponent by $-\varepsilon(1+\delta)$ instead of the expected $-\varepsilon$, as long as 
$\varepsilon < \frac{\gamma}{\delta+1}$ ({\it i.e.} the local exponent remains positive). In the intermediate case $\frac{\gamma}{\delta+1} \leq a \leq \gamma$, the frontier is the union of two line segments. An even more unexpected behaviour then occurs: since its is parallel to the bisector in the neighbourhood of $\sigma=0$, the behaviour of the pointwise exponent is regular ({\it i.e.} it decreases by $\varepsilon$ through $\varepsilon-$differentiation) for $\varepsilon$ small enough. However, when $\varepsilon$ is larger than the ordinate of the point of intersection between the lines $s'\mapsto H(t_0)+s'$ and $s'\mapsto\beta_{t_0}(s')$, namely for $\varepsilon > (\gamma-a)/\delta$, the decrease of the pointwise exponent will be equal to $(\gamma-a)/\delta + (\varepsilon-(\gamma-a)/\delta)(1+\delta)$.
\end{example}

\subsubsection{Uniform almost sure $2$-microlocal frontier of mBm}

Under some assumptions on the function $H$ or its regularity, uniform results hold.
First, in the case where the local regularity of $H$ varies continuously, a direct application of theorem \ref{aslocal} yields the following statement:

\begin{proposition}
Let $X=\left\{X_t;\;t\in\mathbf{R}_{+}\right\}$ be a multifractional Brownian
motion such that the function $t\mapsto\tilde{\beta}(t)$, where
$\tilde{\beta}(t)$ is the local H\"older exponent of $H$ at $t$, is continuous
on some open interval $I$. Then the local H\"older exponent of $X$ satisfies
almost surely
\[
\forall t\in I:\:
\tilde{\alphar}(t)=H(t)\wedge\tilde{\beta}(t).
\]
\end{proposition}

\proof
Using the approximation (\ref{approx}), the deterministic H\"older exponent of
$X$ at $t_0$ can be computed as in proposition \ref{propmbm}
\[
\tilde{\mathbb{\bbalpha}}(t_0)=H(t_0)\wedge\tilde{\beta}(t_0).
\]
The result follows from theorem \ref{aslocal}.
\fin

In the case of a regular mBm, i.e. when the values taken by the function $H$
are smaller than its regularity,
a uniform result for the $2$-microlocal frontier of the process holds as well:

\begin{theorem}\label{thmbmunif}
Let $X=\left\{X_t;\;t\in\mathbf{R}_{+}\right\}$ be a multifractional Brownian
motion such that the function $H$ satisfy, for some open interval $I$,
\[
\forall t\in I;\quad H(t)<\tilde{\beta}(t)
\]
where $\tilde{\beta}(t)$ is the local H\"older exponent of $H$ at $t$.\\
Then, almost surely, the 2-microlocal frontier at all $t_0\in I$ of $X$ in the region
\[
\left\{
\begin{array}{c} 0 < \sigma < 1 + s' \\ -1 < s' < 0 \end{array}
\right.
\]
is equal to the line $\sigma=H(t_0)+s'$.

\noindent In particular, almost surely, for all $t_0 \in I$,
$\alphar(t_0)=\widetilde{\alphar}(t_0)=H(t_0)$.
\end{theorem}

\proof
Under the assumptions of the theorem, for all $t_0\in I$, (\ref{approx}) implies
\begin{equation}\label{mbmminunif}
\forall t_0\in I;\quad\widetilde{\mathbb{\bbalpha}}(t_0)=H(t_0).
\end{equation}

Conversely, for all $t_0\in I$, and all sequence
$\left(\rho_n\right)_{n\in\mathbf{N}}$ converging to $0$, there exist two
sequences $\left(t_n\right)_{n\in\mathbf{N}}$ and
$\left(u_n\right)_{n\in\mathbf{N}}$ such that for all $n\in\mathbf{N}$,
$t_n,u_n\in B(t_0,\rho_n)$ and $|t_n-u_n|=\rho_n$.\\
Then, by (\ref{approx}),
for all $\sigma$ s.t. $H(t_0)+s'<\sigma<\tilde{\beta}(t_0)+s'$, we have
\begin{align*}
\frac{E\left[X_{t_n}-X_{u_n}\right]^2}{|t_n-u_n|^{2\sigma}\rho_n^{-2s'}}
\longrightarrow +\infty
\end{align*}
as $n$ goes to $+\infty$,
which gives
\begin{align*}
\limsup_{\rho\rightarrow 0}\sup_{t,u\in B(t_0,\rho)}
\frac{E\left[X_t-X_u\right]^2}{|t-u|^{2\sigma}\rho^{-2s'}}>0
\end{align*}
and thus
\begin{equation}\label{mbmmajunif}
\forall t_0\in I, \forall s'<0;\quad\mathbb{\bbsigma}_{t_0}(s')\leq H(t_0)+s'.
\end{equation}
The result follows from (\ref{mbmminunif}), (\ref{mbmmajunif}) and corollary
\ref{corfrontunif}.
\fin

\begin{remark}
With global regularity conditions on the function $H$, one can
obtain a uniform analog of Proposition \ref{propmbm}. For
instance, it is not hard to adapt the proofs above to show that if
the inequality
$$ \limsup_{\rho\rightarrow 0}
\sup_{t,u\in B(t_0,\rho)}\frac{|H(t)-H(u)|}{\|t-u\|^{\beta_{t_0}}\rho^{-s'}}
<\infty $$
is verified for all $t_0 \in I$, then, almost surely, for all $t_0 \in I$,
$$\sigmar_{t_0}(s')=(H(t_0)+s') \wedge \beta_{t_0}(s')$$
and
$$\alphar(t_0)=H(t_0) \wedge \beta(t_0)$$
$$\widetilde{\alphar}(t_0)=H(t_0) \wedge \tilde{\beta}(t_0)$$
where $\beta(t_0)$ and $\tilde{\beta}(t_0)$
denote the pointwise and local \ho exponents of $H$ at $t_0$.
\end{remark}

\subsection{Generalized Weierstrass function}

Let us recall the definition of the well-known Weierstrass function (\cite{falc}):
\begin{equation}
W_H(t)=\sum_{j=1}^{\infty}\lambda^{-j H}\sin\lambda^{j}t
\end{equation}
where $\lambda \geq 2$ and $H\in (0,1)$.\\
The \ho regularity of a stochastic version of 
the Weierstrass function has been studied
in \cite{cras}. 
Let $\left(Z_j\right)_{j\in\mathbf{N}}$ be a sequence of
$\mathcal{N}(0,1)$ i.i.d. random variables
and define the {\em generalized Weierstrass function (GW)} as the
following Gaussian process $X=\left\{X_t;t\in\mathbf{R}_{+}\right\}$:
\begin{equation}\label{defWeier}
X_t=W_{H(t)}(t)=\sum_{j=1}^{\infty}Z_j \lambda^{-j H(t)}\sin\lambda^{j}t 
\end{equation}
where $t\mapsto H(t)$ takes values in $(0,1)$.\\
The regularity of this process can be
obtained by the computation of the incremental covariance. 
It is easy to show that $X$ is H\"older-continuous but not
differentiable.

\subsubsection{Bound for the incremental covariance of GW}

\begin{proposition}\label{majcovGW}
Let $X=\left\{X_t;\;t\in\mathbf{R}_{+}\right\}$ be a generalized Weierstrass
function. For all $0<a<b$, there exists positive constants $K=K(a)$ 
and $L=L(a)$ such that
\begin{equation}
\forall t,u\in [a,b];\quad 
E\left[X_t-X_u\right]^2 \leq K |t-u|^{H(t)+H(u)}+L \left(H(t)-H(u)\right)^2.
\end{equation}
\end{proposition}

\proof
From (\ref{defWeier}), one computes
\begin{align}\label{covincGW}
E\left[X_t-X_u\right]^2=&\sum_{j,k} \left(\lambda^{-j H(t)}\sin\lambda^{j}t
-\lambda^{-j H(u)}\sin\lambda^{j}u\right)\nonumber\\
&\times\left(\lambda^{-k H(t)}\sin\lambda^{k}t
-\lambda^{-k H(u)}\sin\lambda^{k}u\right). E\left[Z_j Z_k\right]\nonumber\\
=&\sum_{j} \left(\lambda^{-j H(t)}\sin\lambda^{j}t
-\lambda^{-j H(u)}\sin\lambda^{j}u\right)^2.
\end{align}
Using the decomposition
\[
\lambda^{-j H(t)}\sin\lambda^{j}t -\lambda^{-j H(u)}\sin\lambda^{j}u
=\left(\lambda^{-j H(t)}-\lambda^{-j H(u)}\right)\sin\lambda^{j}t
+\lambda^{-j H(u)}\left(\sin\lambda^{j}t-\sin\lambda^{j}u\right),
\]
we get
\begin{align}\label{decomp}
E\left[X_t-X_u\right]^2 \leq 2\sum_{j}
\left(\lambda^{-j H(t)}-\lambda^{-j H(u)}\right)^2\sin^2\lambda^{j}t
+\;2\sum_{j} \lambda^{-2j H(u)}\left(\sin\lambda^{j}t-\sin\lambda^{j}u\right)^2.
\end{align}
First, let us give an upper bound for the first term of (\ref{decomp}).
\begin{align*}
\sum_{j=1}^{\infty}\left(\lambda^{-j H(t)}-\lambda^{-j H(u)}\right)^2\sin^2\lambda^{j}t
\leq\sum_{j=1}^{\infty}\left(\lambda^{-j H(t)}-\lambda^{-j H(u)}\right)^2.
\end{align*}
By the finite increments theorem, there exists $\tau$ between $H(t)$ and $H(u)$
such that
\begin{align*}
\lambda^{-j H(t)}-\lambda^{-j H(u)}=-j\lambda^{-j\tau} \left(H(t)-H(u)\right)
\log\lambda
\end{align*}
therefore
\begin{align}\label{firsterm}
\sum_{j=1}^{\infty}\left(\lambda^{-j H(t)}-\lambda^{-j H(u)}\right)^2\sin^2\lambda^{j}t
\leq\left(H(t)-H(u)\right)^2\log^2\lambda.\sum_{j=1}^{\infty}j^2\lambda^{-2a.j}.
\end{align}
To deal with the second term of (\ref{decomp}),
for given $t,u\in[a,b]$, we consider the integer $N$ such that
$\lambda^{-(N+1)}\leq |t-u|\leq\lambda^{-N}$.
\begin{equation}
\sum_{j=1}^{\infty} \lambda^{-2j H(u)}\left(\sin\lambda^{j}t-\sin\lambda^{j}u\right)^2
\leq\sum_{j=1}^N \lambda^{-2j H(u)}\left(\sin\lambda^{j}t-\sin\lambda^{j}u\right)^2
+4\sum_{j=N+1}^{\infty} \lambda^{-2j H(u)}.
\end{equation}
Then, using the inequality
\begin{align*}
\left(\sin\lambda^j t-\sin\lambda^j u\right)^2&=
4\sin^2\lambda^j\frac{t-u}{2} \cos^2\lambda^j\frac{t+u}{2}\\
&\leq\lambda^{2j}|t-u|^2\leq\lambda^{2j}\lambda^{-2N},
\end{align*}
we get
\begin{align*}
\sum_{j=1}^N \lambda^{-2j H(u)}\left(\sin\lambda^{j}t-\sin\lambda^{j}u\right)^2
&\leq\lambda^{-2N}\sum_{j=1}^N\lambda^{2j(1-H(u))}\\
&\leq\lambda^{-2N}\lambda^{2(1-H(u))}
\frac{\lambda^{2(N-1)(1-H(u))}}{\lambda^{2(1-H(u))}-1}\\
&\leq\frac{\lambda^{-2NH(u)}}{\lambda^{2(1-H(u))}-1}\\
&\leq\frac{\lambda^{2H(u)}}{\lambda^{2(1-H(u))}-1}.|t-u|^{2H(s)}.
\end{align*}
Moreover, as
\begin{align*}
\sum_{j=N+1}^{\infty} \lambda^{-2j H(u)}&=
\frac{\lambda^{-2(N+1)H(u)}}{1-\lambda^{-2H(u)}}\\
&\leq\frac{|t-s|^{2H(u)}}{1-\lambda^{-2H(u)}},
\end{align*}
we have
\begin{equation}\label{secondterm}
\sum_{j=1}^{\infty} \lambda^{-2j H(u)}\left(\sin\lambda^{j}t-\sin\lambda^{j}u\right)^2
\leq \left(\frac{\lambda^{2H(u)}}{\lambda^{2(1-H(u))}-1}+
\frac{4}{1-\lambda^{-2H(u)}}\right).|t-u|^{2H(u)}.
\end{equation}
As
\[
|t-u|^{2H(u)}-|t-u|^{H(u)+H(t)}=O_{a,b}\left(|t-u|^{H(t)+H(u)}(H(t)-H(u))\right)+O_{a,b}\left(H(t)-H(u)\right)^2,
\]
the result follows from (\ref{decomp}), (\ref{firsterm}) and (\ref{secondterm}).
\fin

To get a upper bound for the $2$-microlocal frontier of the generalized
Weierstrass function, we need the following statement

\begin{proposition}\label{mincovGW}
Let $X=\left\{X_t;\;t\in\mathbf{R}_{+}\right\}$ be a generalized Weierstrass
function. For all $t_0\in\mathbf{R}_{+}$, there exists two sequences
$\left(t_n\right)_{n\in\mathbf{N}}$ and $\left(u_n\right)_{n\in\mathbf{N}}$
converging to $t_0$ and positive constants $k_1$ and $l_1$
such that
\begin{align*}
\forall n\in\mathbf{N};\quad\left(E\left[X_{t_n}-X_{u_n}\right]^2\right)^{\frac{1}{2}}
\geq k_1 |t_n-u_n|^{\frac{H(t_n)+H(u_n)}{2}}
-l_1 \left|H(t_n)-H(u_n)\right|.
\end{align*}
Moreover, if $H$ admits a positive local H\"older exponent at $t_0$,
there exists positive constants $k_2$ and $l_2$ such that
\begin{align*}
\forall n\in\mathbf{N};\quad
\left(E\left[X_{t_n}-X_{u_n}\right]^2\right)^{\frac{1}{2}}
\geq -\left[k_2 \left|H(t_n)-H(u_n)\right| \log|t_n-u_n|
+l_2\right] |t_n-u_n|^{\frac{H(t_n)+H(u_n)}{2}}.
\end{align*}
\end{proposition}

\proof
From (\ref{covincGW}), for all $t,u\in\mathbf{R}_{+}$ and all $n\in\mathbf{N}$,
we have
\[
E\left[X_t-X_u\right]^2\geq\left(\lambda^{-nH(t)}\sin\lambda^n t
-\lambda^{-nH(u)}\sin\lambda^n u\right)^2.
\]
Using the decomposition
\[
\lambda^{-n H(t)}\sin\lambda^{n}t -\lambda^{-n H(u)}\sin\lambda^{n}u
=\left(\lambda^{-n H(t)}-\lambda^{-n H(u)}\right)\sin\lambda^{n}t
+\lambda^{-n H(u)}\left(\sin\lambda^{n}t-\sin\lambda^{n}u\right)
\]
and the triangular inequality, we get
\begin{equation}\label{mincovGWtrig}
E\left[X_t-X_u\right]^2\geq\left(
\left|\lambda^{-n H(t)}-\lambda^{-n H(u)}\right|.
\left|\sin\lambda^{n}t\right|
- \lambda^{-n H(u)}\left|\sin\lambda^{n}t-\sin\lambda^{n}u\right|
\right)^2.
\end{equation}

For all $t_0\in\mathbf{R}_{+}$, there exists a sequence
$\left(t_n\right)_{n\in\mathbf{N}}$ converging to $t_0$,
and such that $\left|\sin\lambda^{n}t_n\right|>\frac{1}{2}$
for all $n\in\mathbf{N}$.
For instance, let us start from a sequence
$\left(\tilde{t}_n=t_0+\frac{\pi}{\lambda^n}\right)_{n\in\mathbf{N}}$
converging to $t_0$, and set, for all $n$
\begin{align*}
t_n=\left\{\begin{array}{l}
\tilde{t}_n\textrm{ if }\left|\sin\lambda^{n}t_n\right|>\frac{1}{2}\\
\tilde{t}_n+\frac{\pi}{2\lambda^n}\textrm{ otherwise}.
\end{array}\right.
\end{align*}

Moreover, for all $t\in\mathbf{R}_{+}$ and all $n\in\mathbf{N}$, there exists
$h_n$ such that $\lambda^{-(n+1)}\leq h_n\leq\lambda^{-n}$ and
$\left|\sin\lambda^{n}(t+h_n)-\sin\lambda^{n}t\right|\geq\frac{1}{10}$.\\
As a consequence, setting $u_n=t_n+h_n$ for all $n$, we get a sequence
$\left(u_n\right)_{n\in\mathbf{N}}$ converging to $t_0$, and such that
\begin{equation*}
\forall n\in\mathbf{N};\;
\left(E\left[X_{t_n}-X_{u_n}\right]^2\right)^{\frac{1}{2}}\geq
\frac{1}{10}\lambda^{-n H(u_n)}
-\left|\lambda^{-n H(t_n)}-\lambda^{-n H(u_n)}\right|
\end{equation*}
and
\begin{equation*}
\forall n\in\mathbf{N};\;
\left(E\left[X_{t_n}-X_{u_n}\right]^2\right)^{\frac{1}{2}}\geq
\frac{1}{2}\left|\lambda^{-n H(t_n)}-\lambda^{-n H(u_n)}\right|
-2\lambda^{-n H(u_n)}.
\end{equation*}
Recall that the local H\"older exponent of $H$ at $t_0$ is
\[
\tilde{\beta}(t_0)=\sup\left\{\beta;\;\limsup_{\rho\rightarrow 0}
\sup_{t,u\in B(t_0,\rho)}\frac{|H(t)-H(u)|}{|t-u|^{\beta}}<+\infty\right\}.
\]
As $\tilde{\beta}(t_0)>0$, we can choose $0<\beta<\tilde{\beta}(t_0)$. This
implies
\begin{align*}
n\left(H(t_n)-H(u_n)\right)&\sim -\left(H(t_n)-H(u_n)\right)\log|t_n-u_n|\\
&\sim -\underbrace{\frac{H(t_n)-H(u_n)}{|t_n-u_n|^{\beta}}}_{\rightarrow 0}
\underbrace{|t_n-u_n|^{\beta}\log|t_n-u_n|}_{\rightarrow 0}\\
&\longrightarrow 0.
\end{align*}
A Taylor expansion gives
\[
\lambda^{-n H(t_n)}-\lambda^{-n H(u_n)}
= \left(H(t_n)-H(u_n)\right) n\log\lambda.\lambda^{-nH(u_n)}
+O\left[\left(H(t_n)-H(u_n)\right)^2 n^2 \lambda^{-nH(u_n)}\right].
\]
Therefore, using $|t_n-u_n|\leq\lambda^{-n}$ and the boundedness of
$n\lambda^{-nH(u_n)}$, there exists $l_1>0$ such that
\begin{equation*}
\forall n\in\mathbf{N};\quad
\left(E\left[X_{t_n}-X_{u_n}\right]^2\right)^{\frac{1}{2}}\geq
\frac{1}{10}|t_n-u_n|^{H(u_n)}
-l_1\;\left|H(t_n)-H(u_n)\right|
\end{equation*}
and there exists $k_2>0$ such that
\begin{align*}
&\forall n\in\mathbf{N};\quad\\
&\left(E\left[X_{t_n}-X_{u_n}\right]^2\right)^{\frac{1}{2}}\geq
k_2\left|H(t_n)-H(u_n)\right|\times
|t_n-u_n|^{H(u_n)}\left(-\log|t_n-u_n|\right)
-2\;|t_n-u_n|^{H(u_n)}.
\end{align*}
We conclude in the same way as in the proof of proposition \ref{majcovGW}.
\fin

\subsubsection{Almost sure $2$-microlocal frontier of GW}

Propositions \ref{majcovGW} and \ref{mincovGW} allow to obtain the
almost sure $2$-microlocal frontier of the generalized Weierstrass function
when $H$ is regular. The situation here is similar to the one of mBm.

\begin{theorem}\label{thGWunif}
Let $X=\left\{X_t;\;t\in\mathbf{R}_{+}\right\}$ be a generalized Weierstrass
function such that the function $H$ satisfy, for some open interval $I$,
\[
\forall t\in I;\quad H(t)<\tilde{\beta}(t)
\]
where $\tilde{\beta}(t)$ is the local H\"older exponent of $H$ at $t$.\\
Then, almost surely, the 2-microlocal frontier at any $t_0\in I$ of $X$ in the region
\[
\left\{
\begin{array}{c} 0 < \sigma < 1+s' \\ -1 < s' < 0 \end{array}
\right.
\]
is equal to the line $\sigma=H(t_0)+s'$.
\end{theorem}

\proof
The proof is similar to the one in the case of mBm. We sketch it below.\\
For each $s'\in (-\infty; 0)$, we introduce
\begin{equation}
\beta_{t_0}(s')=\sup\left\{\beta; \limsup_{\rho\rightarrow 0}
\sup_{t,u\in B(t_0,\rho)}\frac{|H(t)-H(u)|}{\|t-u\|^{\beta}\rho^{-s'}}
<\infty \right\}.
\end{equation}
From (\ref{locfrpoint}), for all $s'\in(-1,0]$, we have
$H(t_0)+s'<\beta_{t_0}(s')$.
Then, for all $\sigma<H(t_0)+s'$, we have
\begin{equation*}
\frac{\|t-u\|^{2H(t_0)}}{\|t-u\|^{2\sigma}\rho^{-2s'}}=
\frac{\|t-u\|^{2(H(t_0)-\sigma)}}{\rho^{-2s'}}\leq
\frac{(2\rho)^{2(H(t_0)-\sigma)}}{\rho^{-2s'}}\rightarrow 0
\end{equation*}
and
\begin{equation*}
\frac{\left[H(t)-H(u)\right]^2}{\|t-u\|^{2\sigma}\rho^{-2s'}}\rightarrow 0.
\end{equation*}
Thus, by proposition \ref{majcovGW}, we have
$\sigma\leq\mathbb{\bbsigma}_{t_0}(s')$. This implies
\begin{equation}\label{GWmajbb}
H(t_0)+s'\leq\mathbb{\bbsigma}_{t_0}(s').
\end{equation}

Conversely,
for all $\sigma$ s.t. $H(t_0)+s'<\sigma<\beta_{t_0}(s')$, the two sequences
$\left(t_n\right)_{n\in\mathbf{N}}$ and
$\left(u_n\right)_{n\in\mathbf{N}}$ given by proposition \ref{mincovGW} can be
chosen such that for all $n\in\mathbf{N}$,
$t_0+\frac{\pi}{\lambda^n}\leq t_n\leq t_0+\frac{3\pi}{2\lambda^n}$
and $\lambda^{-(n+1)}\leq u_n-t_n\leq\lambda^{-n}$.\\
As a consequence, we have
\begin{align*}
\forall n\in\mathbf{N};\quad & t_0\leq t_n\leq u_n \textrm{ and}\\
&\lambda^{-n}\left(\frac{1}{\lambda}+\pi\right)\leq
u_n-t_0\leq\lambda^{-n}\left(1+\frac{3\pi}{2}\right).
\end{align*}
Therefore, setting $\rho_n=u_n-t_0$, we have $t_n,u_n\in B(t_0,\rho_n)$ and
\begin{align*}
\frac{\|t_n-u_n\|^{2H(t_0)}}{\|t_n-u_n\|^{2\sigma}\rho_n^{-2s'}}\geq
\frac{\lambda^{-2H(t_0)(n+1)}}
{\lambda^{-2n\sigma}\lambda^{2ns'}\left(1+\frac{3\pi}{2}\right)^{-2s'}}
=\lambda^{2n(H(t_0)+s'-\sigma)}\rightarrow +\infty.
\end{align*}
As, on the other hand,
\begin{align*}
\frac{\left[H(t_n)-H(u_n)\right]^2}{\|t_n-u_n\|^{2\sigma}\rho_n^{-2s'}}\rightarrow 0
\end{align*}
from proposition \ref{mincovGW}, we get
$\sigma\geq\mathbb{\bbsigma}_{t_0}(s')$.
This implies
\begin{equation}\label{GWminbb}
\mathbb{\bbsigma}_{t_0}(s')\leq H(t_0)+s'.
\end{equation}

From (\ref{GWmajbb}) and (\ref{GWminbb}), we get
\begin{align*}
\forall t_0\in I,\forall s'<0;\quad &\mathbb{\bbsigma}_{t_0}(s')=H(t_0)+s'\\
&\mathbb{\bbalpha}(t_0)=H(t_0).
\end{align*}
Corollary \ref{corfrontunif} then gives the result.
\fin

When $H(t) > \tilde{\beta}(t)$, we are not able to conclude in general but we get pointwise almost sure bounds for the $2$-microlocal frontier. However, it should be possible to obtain a complete almost sure result if one uses the definition of the stochastic Weierstrass function used in \cite{cras} instead of (\ref{defWeier}). The trick consists in summing over a particular set of indices that grows sufficiently fast to infinity, rather than on the whole of $\mathbf{N}$. See \cite{cras} for details.

\begin{proposition}
The $2$-microlocal frontier at any $t_0$ of the generalized Weierstrass function
in the region
\[
\left\{
\begin{array}{c} 0 < \sigma < 1+s' \\ -1 < s' < 0 \end{array}
\right.
\]
is, almost surely, ``above the minimum'' of the line $s'\mapsto H(t_0)+s'$ and the
$2$-microlocal frontier of $H$.
\end{proposition}

\proof We have to distinguish between the following two cases:
\begin{itemize}
\item If $H(t_0)+s'<\beta_{t_0}(s')$, theorem \ref{thGWunif} gives the result.
\item If $\beta_{t_0}(s')<H(t_0)+s'$\\
For all $\sigma<\beta_{t_0}(s')$, we have
\begin{align*}
\frac{\|t-u\|^{2H(t_0)}}{\|t-u\|^{2\sigma}\rho^{-2s'}}\rightarrow 0
\end{align*}
and
\begin{align*}
\frac{\left[H(t)-H(u)\right]^2}{\|t-u\|^{2\sigma}\rho^{-2s'}}\rightarrow 0.
\end{align*}
Then, by proposition \ref{majcovGW} and corollary \ref{cormingauss}, we have
$\sigma\leq\sigmar_{t_0}(s')$ almost surely. This implies almost surely
\begin{equation}
\beta_{t_0}(s')\leq\sigmar_{t_0}(s').
\end{equation}
\end{itemize}
\fin

\begin{remark}
Conversely, for all $\sigma$ s.t. $\beta_{t_0}(s')<\sigma<H(t_0)+s'$, there exist sequences
$(\rho_n)_n$, $(t_n)_n$ and $(u_n)_n$ such that $\forall n; t_n,u_n\in
B(t_0,\rho_n)$ and
\begin{align*}
\frac{\left[H(t_n)-H(u_n)\right]^2}{\|t_n-u_n\|^{2\sigma}\rho_n^{-2s'}}
\rightarrow +\infty.
\end{align*}
Moreover, we have
\begin{align*}
\frac{\|t_n-u_n\|^{2H(t_0)}}{\|t_n-u_n\|^{2\sigma}\rho_n^{-2s'}}\rightarrow 0.
\end{align*}
However, the inequalities given by proposition \ref{mincovGW} are not satisfied by the sequences $(t_n)_n$ and $(u_n)_n$. Thus this cannot be used to get an almost sure upper bound for $\sigmar_{t_0}(s')$. Again, using the definition set in \cite{cras} instead of (\ref{defWeier}) should allow to conclude in general.
\end{remark}

\subsection{Application to Wiener and stable integrals}\label{intstoch}

Let us go back to the example given at the beginning of this work. 
We wish to apply our results to Wiener integrals.
Theorem \ref{thone} applies to such processes and allows to evaluate their almost sure $2$-microlocal frontier at any point.
Before we proceed, we need to set a definition:

\begin{definition}\label{pseudo}
Let $\varphi$ be a deterministic function.
The {\em pseudo-2-microlocal frontier} of $\varphi$ at $t_0$ is the function
$s' \mapsto \Sigma_{t_0}(s')$, defined
for  $s'\in(-\infty;0)$ by
\begin{equation}\label{lafaussedef}
\Sigma_{t_0}(s')=\sup\left\{\sigma;\; \limsup_{\rho\rightarrow 0}
\sup_{t,u\in B(t_0,\rho)}\frac{|\varphi(t)-\varphi(u)|}{\|t-u\|^{\sigma}\rho^{-s'}}
<\infty \right\}.
\end{equation}
\end{definition}

In other words, the pseudo-frontier of $\varphi$ is obtained by using definition \ref{def2mlsimple} for all $s'<0$ and all $\sigma$, instead of using definition \ref{def2ml} when $(s',\sigma) \notin D$. In general, the ``true'' 2-microlocal frontier and the pseudo one do not coincide. However, one can show that if either one of these frontiers passes through $D$, then also does the other one, and in this case they coincide in $D$ (see 
\cite{theseAE}, Proposition 3.14 p.43 and also Proposition 3.15, p.44 for a more general result).
A  simple example where the frontiers differ everywhere is provided by the function $\varphi(t) = t + |t|^\gamma$, where $1<\gamma<2$.
The frontier of $\varphi$ at 0 is the line $\sigma(s') = \gamma + s'$ (the regular part is ignored in definition \ref{def2ml}), while its pseudo-frontier is given by $\Sigma(s') = 1+s'$. Note that an even more radical difference is observed by simply taking $\varphi(t) = t$. In contrast, for a chirp, both frontiers coincide, whether they intersect $D$ or not.

\begin{theorem}\label{WI}
Let $X$ be the stochastic process defined by
\begin{equation*}
X_t=\int_0^t \eta(u).dB_u + \psi(t),
\end{equation*}
where $\eta$ and $\psi$ are $L^2$ deterministic functions and $B$ is standard Brownian motion. 
Let $s' \mapsto \beta_{t_0}(s')$ (resp.
$s'\mapsto\gamma_{t_0}(s')$) denote the pseudo-$2$-microlocal frontier of the function $\varphi: t\mapsto\int_0^t \eta^2$ (resp. $\psi$) at $t_0$. \\
Then, the $2$-microlocal frontier of $X$ at $t_0$ in the region
defined by conditions (\ref{detsp}) is almost surely equal to
$\frac{1}{2}\beta_{t_0}(2s')\wedge\gamma_{t_0}(s')$, 
provided this function intersects the considered region.
\end{theorem}

\proof Write
\[
\forall s,t\in\textbf{R}_{+};\quad
X_t-X_s=\int_s^t \eta(u).dB_u + \psi(t)-\psi(s).
\]
By definition of the Wiener integral,
\begin{eqnarray*}
\forall s < t;\quad 
E\left[X_t-X_s\right]^2&=&\left|\int_s^t \eta^2(u).du\right|
+\left(\psi(t)-\psi(s)\right)^2\\
&=&|\varphi(t)-\varphi(s)|+\left(\psi(t)-\psi(s)\right)^2.
\end{eqnarray*}
For all $\rho>0$, all $\sigma>0$ and all $s'<0$,
\[
\forall s,t\in B(t_0,\rho);\;
\frac{E\left[X_t-X_s\right]^2}{|t-s|^{2\sigma}\rho^{-2s'}}
=\frac{|\varphi(t)-\varphi(s)|}{|t-s|^{2\sigma}\rho^{-2s'}}
+\left(\frac{|\psi(t)-\psi(s)|}{|t-s|^{\sigma}\rho^{-s'}}\right)^2.
\]
The deterministic $2$-microlocal frontier of $X$ at $t_0$ is
therefore
\begin{eqnarray*}
\mathbb{\bbsigma}_{t_0}(s')=\frac{1}{2}\beta_{t_0}(2s')\wedge\gamma_{t_0}(s').
\end{eqnarray*}

Since $X$ is a Gaussian process, the result follows from Theorem \ref{thone}. \fin

\begin{remark}
Assume $\psi \equiv 0$ so that the frontier of $X$ is equal to $\frac{1}{2}\beta_{t_0}(2s')$.
Using the characterizations of the pointwise and local exponents in terms of the frontier, one can easily see that $\widetilde{\alphar}(t_0) = \widetilde{\alpha}_\varphi(t_0)/2$ and $\alphar(t_0) = \alpha_\varphi(t_0)/2$ almost surely, at least when the pseudo-frontier coincides with the true frontier.
\end{remark}

\begin{example} \label{ExWi}
Following the introduction section, consider a Wiener integral $X$ whose kernel $\eta$ is
equal to $\sqrt{|t-t_0|^\gamma |\sin(|t-t_0|^{-\delta})|}$ in
a neighbourhood of $t_0$, where $\gamma>0, \delta>0$ (we take $\psi\equiv0$).
We may apply theorem \ref{WI} to get the almost
sure 2-microlocal frontier of $X$ at $t_0$. 
The pseudo-2-microlocal frontier of $\varphi$ (as defined in theorem \ref{WI}) 
is equal to its plain 2-microlocal frontier, and is given by 
$\beta_{t_0}(s')=\frac{1} {\delta+1} s' +
\frac{\gamma}{\delta+1}+1$ (since $\varphi$ is a primitive of a chirp, its 
frontier is the one of the chirp plus one). Thus, the frontier of $X$ is equal to 
$\frac{1} {\delta+1} s' +
\frac{\gamma}{2\delta+2}+\frac{1}{2}$. It is then straightforward to check
the values of the pointwise and local exponents announced in the introduction. 
In addition, one can see that differentiation of $X$ of order $\varepsilon$ will
decrease its pointwise exponent by $\varepsilon(1+\delta)$, as long as 
the local exponent remains positive. 

In this case, one could have considered the true 2-microlocal frontier of $\varphi$ as, for a chirp, it coincides with the pseudo-frontier. 
But the use of the true $2$-microlocal frontier of $\varphi$, instead of its pseudo-frontier, can lead to wrong prediction. 
For instance, assume that we replace $\eta(t)=\sqrt{|t-t_0|^\gamma |\sin(|t-t_0|^{-\delta})|}$ with $\mu(t) = \sqrt{1+|t-t_0|^\gamma |\sin(|t-t_0|^{-\delta})|}$ and consider the Wiener integral $Y_t$ whose kernel is $\mu$. Using the 2-microlocal frontier
would lead to the prediction that pointwise and local exponents of $Y$ should be the same as the ones of $X$ (since the frontier of $\mu$ and $\varphi$ are the same). 
In contrast, the use of the pseudo-frontier of $\mu$ yields the correct result that $Y$ has both local and pointwise exponents equal to 1/2.
\end{example}

\begin{remark}
In the above example, when $\gamma,\delta$ tend to 0, the 2-microlocal frontier of $X$ tends to the
one of Brownian motion, and of course, both the exponents tend to $1/2$. This is
consistent with the fact that $X$ tends to Brownian motion. This ``continuity''
property is not true in the deterministic frame: when $\gamma,\delta$ tend to 0,
$\varphi$ tends to a linear function, whose frontier is not the limit of the frontier
of $\varphi$, which is equal to $\beta(s')=s'+1$. 
On the contrary to the ``true" frontier, this property holds for the pseudo-frontier.
\end{remark}

Using proposition \ref{propmin} instead of theorem \ref{WI}, one may compute
in the same manner as above a lower bound for the 2-microlocal frontier of stable integrals.
Such integrals are obtained by replacing the Wiener measure by a stable measure
in the definition of $X$ (see \cite{ST} for an account on stable integrals).
More precisely, let $X$ be the stochastic process defined by
\begin{equation*}
X_t=\int_0^t \eta(u).dM_u ,
\end{equation*}
where $\eta$ is an $L^\alpha$ deterministic function and $M$ is an $\alpha-$stable
random measure ($\alpha \in (0,2)$). Then, provided $\alpha \neq 1$,
or, if $\alpha=1$ assuming $M$ is symmetric,

\begin{eqnarray*}
\forall s < t, \forall p \in (0,\alpha);\quad 
E\left|X_t-X_s\right|^p&=&C(\alpha,p)\left(\int_s^t |\eta(u)|^\alpha.du\right)^{\frac{p}{\alpha}},\\
\end{eqnarray*}
where $C(\alpha,p)$ is a constant depending only on $\alpha$ and $p$ (see \cite{ST}).

Reasoning as in the Gaussian case, one
may obtain a lower bound on the frontier of $X$ in terms $\beta_{t_0}$, 
the pseudo-$2$-microlocal frontier of the function 
$\varphi: t\mapsto\int_0^t |\eta|^\alpha$. Indeed,
\begin{eqnarray*}
\forall s,t \in B(t_0,\rho), \forall p \in (0,\alpha);\quad 
\frac{E\left|X_t-X_s\right|^p}{|t-s|^{\mu+1}\rho^{-\nu}}&=&C(\alpha,p)\left(\frac{|\varphi(t) -\varphi(s)|}{|t-s|^{(\mu+1)\frac{\alpha}{p}}\rho^{-\nu\frac{\alpha}{p}}}\right)^{\frac{p}{\alpha}}.\\
\end{eqnarray*}
By definition of $\beta_{t_0}$, the right-hand side in the equality above will be finite
when 
\begin{eqnarray*}
(\mu+1)\frac{\alpha}{p} < \beta_{t_0}(\nu\frac{\alpha}{p}).
\end{eqnarray*}
Using proposition \ref{propmin}, this entails that
\begin{eqnarray*}
\sigmar_{t_0}(\frac{\nu}{p})\geq \frac{1}{\alpha}\beta_{t_0}(\frac{\nu\alpha}{p})- \frac{1}{p}.
\end{eqnarray*}
Now set $s'=\frac{\nu}{p}$ to get:
\begin{eqnarray*}
\sigmar_{t_0}(s')\geq \frac{1}{\alpha}\beta_{t_0}(\alpha s')- \frac{1}{p}.
\end{eqnarray*}
Since this is true for all $p \in (0,\alpha)$, we get finally:
\begin{eqnarray*}
\sigmar_{t_0}(s')\geq \frac{1}{\alpha}\beta_{t_0}(\alpha s')- \frac{1}{\alpha}.
\end{eqnarray*}
Although this result is less precise than in the Gaussian case, a lower bound for the local regularity is a very interesting result for sample paths of a stochastic process.

\begin{example} \label{ExStable}

Consider a stable integral $X$ whose kernel is
equal to $\left(|t-t_0|^\gamma |\sin(|t-t_0|^{-\delta})|\right)^{\frac{1}{\alpha}}$ in
a neighbourhood of $t_0$, where $\gamma>0, \delta>0$.
We may apply the result above to get that
the frontier of $X$ is almost surely not smaller than 
$\frac{1} {\delta+1} s' +
\frac{\gamma}{\alpha(\delta+1)}$. This is the frontier of a chirp
$t \mapsto |t-t_0|^{\frac{\gamma}{\alpha}}\sin(|t-t_0|^{-\delta})$.

We also mention the following curiosity: when $\alpha=1$ and the
skewness function $\lambda$ of $M$ is not zero, one has
\begin{eqnarray*}
E\left|X_t-X_s\right|^p&=&C(p)\left(\int_s^t |\eta(u)|.du + \int_s^t |\eta(u)\lambda(u) \log(|\eta(u)|)|.du\right)^p,\\
\end{eqnarray*}
for all $p \in (0,1)$ In this case, one can obtain the same kind of 2-microlocal behaviour
as above even with a smooth function $\eta$ by putting all the irregularity in 
the skewness function $\lambda$.

\end{example}

\section{Proof of intermediate results}\label{boundsect}

\subsection{Proof of Proposition \ref{propunifminsigma}}
Let us suppose first that $s'$ is fixed and that
the function $\varsigma$ is constant equal to $\sigma$.\\
Let $\epsilon>0$ such that $\tilde{\sigma}=\sigma-\epsilon>0$.
Let us take $\rho=2^{-n}$ for $n\geq n_0=-\log_2\delta$ and set
$D_n^m(t_0)=\left\{t_0+k.2^{-(m+n)};k\in\left\{0, \pm 1,\dots,\pm (2^m-1)\right\}^N\right\}$.
Let us consider the event
\[
\Omega^m_n=\left\{\max_{\scriptstyle i\in\left\{0,\dots,2^{m+n}\right\}^N \atop
{\scriptstyle k,l\in\left\{0,\dots,\pm 2^m\right\}^N
\atop \scriptstyle \|k-l\|=1}}
\left|X_{(i+k).2^{-(m+n)}}-X_{(i+l).2^{-(m+n)}}\right|>
2^{-\tilde{\sigma}(m+n)}2^{s'n}\right\}.
\]
For all $p\in\mathbf{N}^{*}$, we have
\begin{align*}
P\left\{\Omega^m_n\right\}
&\leq\sum_{\scriptstyle i\in\left\{0,\dots,2^{m+n}\right\}^N \atop
{\scriptstyle k,l\in\left\{0,\dots,\pm 2^m\right\}^N
\atop \scriptstyle \|k-l\|=1}} P\left\{
\left|X_{(i+k).2^{-(m+n)}}-X_{(i+l).2^{-(m+n)}}\right|>
2^{-\tilde{\sigma}(m+n)}2^{s'n}\right\}\\
&\leq\sum_{\scriptstyle i\in\left\{0,\dots,2^{m+n}\right\}^N \atop
{\scriptstyle k,l\in\left\{0,\dots,\pm 2^m\right\}^N
\atop \scriptstyle \|k-l\|=1}}\frac{E\left[X_{(i+k).2^{-(m+n)}}-X_{(i+l).2^{-(m+n)}}\right]^{2p}}
{2^{-2p\tilde{\sigma}(m+n)}2^{2ps'n}}\\
&\leq 2N\;C\;\lambda_p\;\underbrace{\#\left\{0,\dots,\pm 2^m\right\}^N}_{\left(1+2^{m+1}\right)^N}
\;\underbrace{\#\left\{0,\dots,2^{m+n}\right\}^N}_{\left(1+2^{m+n}\right)^N}
2^{-2p\epsilon (m+n)},
\end{align*}
where $\lambda_p$ is the positive constant such that for all centered
Gaussian random variable $Y$ and all $p\in\mathbf{N}^{*}$, we have
$E\left[Y^{2p}\right]=\lambda_p \left(E\left[Y^2\right]\right)^p$.\\
Then,
\begin{equation*}
P\left\{\Omega^m_n\right\}
\leq (2N.2^{3N})\;C\;\lambda_p\; 2^{(2N-2p\epsilon)m} 2^{(N-2p\epsilon) n}.
\end{equation*}
Choosing $p\in\mathbf{N}$ such that $2N-2p\epsilon<0$,
we deduce
\begin{align*}
P\left\{\exists m;\Omega^m_n\right\}=P\left\{\bigcup_m \Omega^m_n\right\}
&\leq\sum_m P\left\{\Omega^m_n\right\}\\
&\leq\frac{(2N.2^{3N})\;C\;\lambda_p\;2^{(N-2p\epsilon) n}}{1-2^{2N-2p\epsilon}}.
\end{align*}
The Borel-Cantelli lemma implies existence of a random variable
$n^{*}\geq n_0$ such that, almost surely,
\begin{equation}\label{BCunif}
\forall n\geq n^{*}, \forall m\in\mathbf{N};\quad
\max_{\scriptstyle i\in\left\{0,\dots,2^{m+n}\right\}^N \atop
{\scriptstyle k,l\in\left\{0,\dots,\pm 2^m\right\}^N
\atop \scriptstyle \|k-l\|=1}}
\left|X_{(i+k).2^{-(m+n)}}-X_{(i+l).2^{-(m+n)}}\right|
\leq 2^{-\tilde{\sigma}(m+n)}2^{s'n}.
\end{equation}
Therefore, setting $E_r=\left\{i.2^{-r};\; i\in\left\{0,\dots,2^{r}\right\}^N\right\}\subset [0,1]$,
we show that for all $n\geq n^{*}$ and all $m\in\mathbf{N}$,
\begin{align}\label{recunif}
\forall q>m,\; \forall t_0\in E_{q+n};\;
\forall t,u\in D_n^q(t_0)\textrm{ s.t. } \|t-u\|<2^{-(m+n)};\nonumber\\
|X_t-X_u|&\leq
2\left(\sum_{j=m+1}^{q}2^{-\tilde{\sigma}(j+n)}\right) 2^{s'n}\\
&\leq\frac{2.2^{-\tilde{\sigma}(m+n+1)}}{1-2^{-\tilde{\sigma}}}
2^{s' n}. \nonumber
\end{align}
To prove (\ref{recunif}), we proceed by induction:
\begin{itemize}
\item for $q=m+1$, for all $t_0\in E_{m+n+1}$, the conditions $t,u\in
D_n^{m+1}(t_0)$ and $\|t-u\|<2^{-(m+n)}$ impose on $t$ and $u$ to be
neighbors in $D_n^{m+1}(t_0)$. 
Therefore (\ref{recunif}) follows from (\ref{BCunif}).

\item assume that the property is valid for an integer $M>m$,
then take $t_0\in E_{M+n+1}$, and $t,u\in D_n^{M+1}(t_0)$ such that
$|t-u|<2^{-(m+n)}$.
There exists $\tilde{t_0}\in E_{M+n}$ such that
$|t_0-\tilde{t_0}|\leq 2^{-(M+n+1)}$.
As $\tilde{t_0}$ can be chosen such that the following strict inequality holds
\begin{align*}
\|u-\tilde{t_0}\| &< \|u-t_0\| + \|t_0-\tilde{t_0}\|\\
&< 2^{-n}-2^{-(M+n+1)}+2^{-(M+n+1)},
\end{align*}
there exists $\tilde{u}\in D_n^{M}(\tilde{t_0})$ such that
$|u-\tilde{u}|\leq 2^{-(M+n+1)}$.
In the same way, we get $\tilde{t}\in D_n^{M}(\tilde{t_0})$ such that
$|t-\tilde{t}|\leq 2^{-(M+n+1)}$.
Moreover, $\tilde{t}$ and $\tilde{u}$ can be chosen such that
$|\tilde{t}-\tilde{u}|\leq |t-u| < 2^{-(m+n)}$.
Then, by the triangular inequality
\begin{equation*}
|X_t-X_u|\leq |X_t-X_{\tilde{t}}|+|X_{\tilde{t}}-X_{\tilde{u}}|
+|X_{\tilde{u}}-X_u|
\end{equation*}
and the fact that $t,\tilde{t},\tilde{u},u$ belong to $D_n^{M+1}(\tilde{t_0})$,
(\ref{BCunif}) gives
\begin{equation*}
|X_t-X_u|\leq 2.2^{-\tilde{\sigma}(M+n+1)} 2^{s' n} +
|X_{\tilde{t}}-X_{\tilde{u}}|.
\end{equation*}
\end{itemize}
Property (\ref{recunif}) follows.

Let us take $t_0\in\bigcup_q E_{q+n}$ and $t,u\in\bigcup_q D_n^q(t_0)$.
There exists $m>0$ such that $2^{-(m+n+1)}\leq\|t-u\| < 2^{-(m+n)}$.
Then property (\ref{recunif}) applied to $m$, $t_0$, $t$ and $u$ gives
$|X_t-X_u|\leq\frac{2}{1-2^{-\tilde{\sigma}}}
\|t-u\|^{\tilde{\sigma}}.2^{s'n}$.\\
Using the continuity of $X$, we get
\begin{align*}
\forall t_0\in [0,1];\;
\forall t,u\in B(t_0,2^{-n});\quad |X_t-X_u|\leq
\frac{2}{1-2^{-\tilde{\sigma}}} \|t-u\|^{\tilde{\sigma}}
2^{s'n}.
\end{align*}
Hence, almost surely, for all $\rho\in (0,2^{-n^{*}})$, there exists
$n>n^{*}$ such that $2^{-(n+1)}\leq\rho\leq 2^{-n}$ and
\begin{align}\label{minfrontconst}
\forall t,u\in B(t_0,\rho); \quad 
|X_t-X_u|&\leq\frac{2}{1-2^{-\tilde{\sigma}}}
\|t-u\|^{\tilde{\sigma}} 2^{s'n}\nonumber\\
&\leq 2^{-s'}\frac{2}{1-2^{-\tilde{\sigma}}}
\|t-u\|^{\tilde{\sigma}} \rho^{-s'}.
\end{align}

\vsp
In the general case where $\varsigma$ is not constant,
for all $a,b\in\mathbf{Q}_{+}^N$ with $a\prec b$
let us consider
$\sigma=\inf_{u\in [a,b]}\varsigma_{u}(s')-\epsilon$ with $\epsilon>0$.
By (\ref{minfrontconst}), there exists a set $\Omega^{*}\subset\Omega$
such that $P\left\{\Omega^{*}\right\}=1$ and for all $\omega\in\Omega^{*}$,
\[
\forall a,b\in\mathbf{Q}_{+}^N,\forall\epsilon\in\mathbf{Q}_{+},
\forall t_0\in\overbrace{[a,b]}^{\circ};\quad
\sigmar(t_0)\geq\inf_{u\in [a,b]}\varsigma_{u}(s')-\epsilon.
\]
Therefore, taking two sequences $\left(a_n\right)_{n\in\mathbf{N}}$ and
$\left(b_n\right)_{n\in\mathbf{N}}$ such that $\forall
n\in\mathbf{N};\;a_n<t_0<b_n$ and converging to $t_0$,
we have for all $\omega\in\Omega^{*}$
\begin{equation}
\forall t_0\in\mathbf{R}^N_{+};\quad
\sigmar_{t_0}(s')\geq\liminf_{u\rightarrow t_0}\varsigma_{u}(s').
\end{equation}

\subsection{Proof of Proposition \ref{propmajunif}}
First of all, a classical proof allows to show that lemma \ref{majlem} implies
\[
\forall \tilde{\sigma}>\sigma;\quad
P\left\{\forall t_0\in\mathbf{Q}^N_{+}; \sigmar_{t_0}(s')<\tilde{\sigma}\right\}=1.
\]
To extend this result for all $t_0\in\mathbf{R}^N_{+}$,
let us consider $t_0\in\mathbf{R}^N_{+}-\mathbf{Q}^N_{+}$, and a sequence
$(x^{(m)})_{m\in\mathbf{N}}$ in $\mathbf{Q}^N_{+}$ such that $x^{(m)}\rightarrow
t_0$. We have, almost surely, for all $m$, $\sigmar_{x^{(m)}}(s')<\tilde{\sigma}$.
Then there exists sequences $(\rho^{(m)}_n)_{n\in\mathbf{N}}$,
$(t^{(m)}_n)_{n\in\mathbf{N}}$ and $(u^{(m)}_n)_{n\in\mathbf{N}}$ such that
for all $n\in\mathbf{N}$, $t^{(m)}_n,u^{(m)}_n\in B(x^{(m)},\rho^{(m)}_n)$ and
\begin{equation}\label{liminfty}
\lim_{n\rightarrow\infty}\frac{|X_{t^{(m)}_n}-X_{u^{(m)}_n}|}
{\|t^{(m)}_n-u^{(m)}_n\|^{\tilde{\sigma}}(\rho^{(m)}_n)^{-s'}}=+\infty.
\end{equation}
From these $m$ sequences, we build $3$ sequences $(t_n)$, $(u_n)$ and $(\rho_n)$
such that $t_n\rightarrow t_0$, $u_n\rightarrow t_0$,
$t_n,u_n \in B(t_0,\rho_n)$ and
\[
\lim_{n\rightarrow\infty}\frac{|X_{t_n}-X_{u_n}|}
{\|t_n-u_n\|^{\tilde{\sigma}}(\rho_n)^{-s'}}=+\infty.
\]
For all $n$ and $m$, let us write
\begin{equation}\label{decompGWsuite}
\left\{\begin{array}{l}
t^{(m)}_n-t_0=t^{(m)}_n-x^{(m)}+x^{(m)}-t_0\\
u^{(m)}_n-t_0=u^{(m)}_n-x^{(m)}+x^{(m)}-t_0.
\end{array}\right.
\end{equation}
Let us fix $\epsilon>0$. There exists $N>0$ such that
\begin{equation}\label{suitexm}
\forall m\geq N; \quad |x^{(m)}-t_0|<\epsilon.
\end{equation}
As for all $m$, $|t^{(m)}_n-x^{(m)}|\rightarrow 0$ when $n\rightarrow\infty$,
there exists a subsequence $(t^{(m)}_{p_n})_n$ such that
\[
\forall n;\quad
|t^{(m)}_{p_n}-x^{(m)}|<|t^{(m-1)}_n-x^{(m-1)}|
\]
therefore, we can suppose that the sequences $(t^{(m)}_n)_n$ satisfy the
properties
\begin{equation}\label{decs}
\forall m,n;\quad
|t^{(m)}_n-x^{(m)}|<|t^{(m-1)}_n-x^{(m-1)}|.
\end{equation}
Moreover, in the same way, we can suppose
\begin{equation}\label{dect}
\forall m,n;\quad
|u^{(m)}_n-x^{(m)}|<|u^{(m-1)}_n-x^{(m-1)}|.
\end{equation}
Then, using the fact that there exists $N'\in\mathbf{N}$ such that
\[
\forall n\geq N';\quad \left\{
\begin{array}{l}|t^{(1)}_n-x^{(1)}|<\epsilon\\
|u^{(1)}_n-x^{(1)}|<\epsilon
\end{array}\right.
\]
(\ref{decs}) and (\ref{dect}) imply
\begin{equation}\label{term2}
\forall n\geq N';\quad \left\{
\begin{array}{l}|t^{(n)}_n-x^{(n)}|<\epsilon\\
|u^{(n)}_n-x^{(n)}|<\epsilon.
\end{array}\right.
\end{equation}
Therefore (\ref{decompGWsuite}), (\ref{suitexm}) and (\ref{term2}) lead to
\begin{equation}\label{limt0}
\lim t^{(n)}_n=\lim u^{(n)}_n=t_0.
\end{equation}

In the same way as previously, by (\ref{liminfty}), we can suppose that the
sequences $(t^{(m)}_n)_n$, $(u^{(m)}_n)_n$ and $(\rho^{(m)}_n)_n$ satisfy
\begin{equation}\label{croissinfty}
\forall m,n;\quad
\frac{|X_{t^{(m)}_n}-X_{u^{(m)}_n}|}
{\|t^{(m)}_n-u^{(m)}_n\|^{\tilde{\sigma}}(\rho^{(m)}_n)^{-s'}}
>\frac{|X_{t^{(m-1)}_n}-X_{u^{(m-1)}_n}|}
{\|t^{(m-1)}_n-u^{(m-1)}_n\|^{\tilde{\sigma}}(\rho^{(m-1)}_n)^{-s'}}.
\end{equation}
Then, for all $M>0$, there exists $N''\in\mathbf{N}$ such that
\[
\forall n\geq N'';\quad
\frac{|X_{t^{(1)}_n}-X_{u^{(1)}_n}|}
{\|t^{(1)}_n-u^{(1)}_n\|^{\tilde{\sigma}}(\rho^{(1)}_n)^{-s'}}>M
\]
which leads to, using (\ref{croissinfty})
\begin{equation}
\forall n\geq N'';\quad
\frac{|X_{t^{(n)}_n}-X_{u^{(n)}_n}|}
{\|t^{(n)}_n-u^{(n)}_n\|^{\tilde{\sigma}}(\rho^{(n)}_n)^{-s'}}>M.
\end{equation}
We have shown
\begin{equation}\label{liminftyn}
\lim_{n\rightarrow\infty}\frac{|X_{t^{(n)}_n}-X_{u^{(n)}_n}|}
{\|t^{(n)}_n-u^{(n)}_n\|^{\tilde{\sigma}}(\rho^{(n)}_n)^{-s'}}=+\infty.
\end{equation}

Therefore (\ref{limt0}) and (\ref{liminftyn}) imply
$\sigmar_{t_0}(s')\leq\tilde{\sigma}$. Thus we can state that, almost surely,
\[
\forall t_0\in\mathbf{R}^N_{+};\quad \sigmar_{t_0}(s')\leq\sigma.
\]

\bibliographystyle{plain}
\bibliography{style}

\end{document}